\documentclass[10pt,journal]{IEEEtran}

\usepackage{amsfonts}
\usepackage{extarrows}
\usepackage{enumerate}
\usepackage{amsmath}
\usepackage{algorithm,algorithmic}
\usepackage{amssymb}
\usepackage{subcaption}

\usepackage[colorlinks=true,linkcolor=blue,citecolor=blue]{hyperref}
\usepackage{graphicx}
\usepackage{color}
\usepackage[square,comma,sort&compress,numbers]{natbib}

\allowdisplaybreaks[4]

\begin{document}

% \title{Learning-Based Adaptive Observation and Approximate Optimal Control for Uncertain Systems}

% \title{Adaptive Observation Based Reinforcement Learning for Approximate Optimal Control of Uncertain Systems}

\title{Adaptive Observation-Based Efficient Reinforcement Learning for Uncertain Systems}

 \author{Maopeng Ran and Lihua Xie, \emph{Fellow, IEEE}
 \thanks{This work was supported by the Delta-NTU Corporate Lab through the National Research Foundation Corporate Lab@University Scheme, and the Projects of Major International (Regional) Joint Research Program NSFC (Grant No. 61720106011).}
 \thanks{M. Ran and L. Xie (corresponding author) are with the School of Electrical and Electronic Engineering, Nanyang Technological University, Singapore 639798 (email: mpran@ntu.edu.sg; elhxie@ntu.edu.sg).}}

% \markboth{IEEE TRANSACTIONS ON AUTOMATIC CONTROL}
% \markboth{}
%{Maopeng Ran \MakeLowercase{\textit{et al.}}: ADRC for time-delay systems}
\maketitle

\begin{abstract}

This paper develops an adaptive observation-based efficient reinforcement learning (RL) approach for systems with uncertain drift dynamics.  A novel concurrent learning adaptive extended observer (CL-AEO) is first designed to jointly estimate the system state and parameter. This observer has a two-time-scale structure and doesn't require any additional numerical techniques to calculate the state derivative information.  The idea of concurrent learning (CL) is leveraged to use the recorded data, which leads to a relaxed verifiable excitation condition for the convergence of parameter estimation. Based on the estimated state and parameter provided by the CL-AEO, a simulation of experience based RL scheme is developed to online approximate the optimal control policy. Rigorous theoretical analysis is given to show that the practical convergence of the system state to the origin and the developed policy to the ideal optimal policy can be achieved without the persistence of excitation (PE) condition. Finally, the effectiveness and superiority of the developed methodology are demonstrated via comparative simulations.

\end{abstract}


\IEEEpeerreviewmaketitle

\section{Introduction}

% (RL) -> problem -> model-based RL -> Open Problems -> Our contribution

Reinforcement learning (RL), inspired by learning mechanisms observed in naturally occurring systems (e.g., animals and social groups \cite{Lewis-2013}), is concerned with how agents or actors ought to take optimal actions in an environment to maximize the notion of cumulative reward \cite{Kae-1996}. In the last several decades, RL has been adopted in control theory and has had an increasing success in finding adaptive optimal policies for  dynamic control systems \cite{Lewis-2018,Jiang-2013}.

Early efforts for the implementation of RL algorithms in control society can be dated back to 1980s \cite{Werbos-1989,Werbos-1992}, in which RL was employed to solve the optimal regulation problem for discrete-time systems. Extending RL to continuous-time systems was first visited in \cite{Doya-2000}. After that, considerable RL solutions were developed for  both  discrete-time and continuous-time systems. In \cite{Lewis-2009}, an RL algorithm which solves the algebraic Riccati equation corresponding to
the LQR problem was proposed without requiring the knowledge of the system state matrix. In \cite{Jiang-2012}, the approach in \cite{Lewis-2009} was extended to continuous-time linear systems with completely unknown system dynamics. In \cite{Lewis-2010}, an RL algorithm was developed to solve the infinite horizon optimal control problem for nonlinear systems with known dynamics. The RL algorithm provides an online approximated solution to the Hamilton–Jacobi–Bellman (HJB) equation. The infinite horizon optimal control problem was further investigated for affine and nonaffine unknown nonlinear systems in \cite{Luo-2014} and \cite{Jiang-2014}, respectively. In recent years, RL algorithms were developed for more complex systems such as polynomial systems \cite{Jiang-2015} and nonstrict-feedback systems \cite{Bai-2020}, and more complex tasks such as multi-agent formation \cite{Wen-2017} and fault-tolerant control \cite{Ma-2019}.

In RL-based control, a neural network (NN) based actor-critic structure is generally employed to online approximate the ideal optimal control policy \cite{Wang-2017}. Therefore, unlike traditional adaptive controllers, the ideal weights of the NN must be exactly learned. This puts a significant challenge for the implementation of RL algorithms \cite{Dep-2018}. A common way to handle this challenge is to assume a persistence of excitation (PE) condition. The PE condition guarantees that the system state explores sufficient points in the state space to generate an “ideal” approximation over the entire domain of operation. However, this condition is difficult to be theoretically guaranteed and cannot be monitored online. In many aforementioned RL works \cite{Jiang-2012,Lewis-2010,Jiang-2014,Jiang-2015,Ma-2019}, to fulfill the PE condition,  carefully selected probing signals are injected into the system, which will inevitably cause undesirable oscillations. Due to this issue, data-driven techniques such as experience replay \cite{Lewis-2014,Xue-2020,Yang-2019} were leveraged to improve data efficiency in online approximate optimal control by reusing the recorded data, and consequently to relax the requirement of the PE condition. However, as pointed out in \cite{Kama-2016a}, since the data can only be recorded along the system trajectory, the system under the designed experience replay-based RL controller still needs to provide sufficient excitation for learning. For example,  probing signals were used in the numerical examples in \cite{Lewis-2014}.

Note that in RL, the NN weights are updated using Bellman error (BE) as a performance metric. If the system dynamics is known, the BE can be evaluated at any desired point in the state space, rather than only along the system trajectory. In this case,  sufficient exploration can be guaranteed by appropriately selecting the points to cover the domain of operation. This idea to improve data efficiency is interpreted as simulation of experience and falls into the so-called model-based RL \cite{Kama-2016a,Kama-2016b,Dep-2020,Dep-2020b}. The model-based RL is capable of relaxing the PE condition and removing the requirement of the probing signal. In model-based RL, one of the main tasks is to obtain the model information before or along with the learning process. In \cite{Kama-2016a,Kama-2016b,Dep-2020}, the system drift dynamics was online learned with the knowledge of full state and state derivative. If exact state derivative information is not available, additional numerical smoothing techniques are needed, which will introduce smoothing errors and increase processing and storage burden. In \cite{Dep-2020b}, an integral concurrent learning (ICL) estimator with full state feedback was proposed to learn the system drift dynamics without the requirement of the state derivative information. However, the ICL requires numerical techniques to evaluate the integrals, which will inevitably cause accumulated errors.

Based on the above discussions and considerations, in this paper, we propose an adaptive observation-based approach to enable efficient model-based RL without the state derivative information or integral calculation. Our approach is inspired from the communities of adaptive observers \cite{Farza-2018,Tyu-2013,Ibr-2018}, concurrent learning \cite{Chow-2013}, and model-based RL \cite{Kama-2016a,Kama-2016b,Dep-2020,Dep-2020b}. The main contributions of this paper are twofold:

\begin{enumerate}[1)]
  \item A concurrent learning adaptive extended observer (CL-AEO) is proposed for joint state and parameter estimation. This observer falls into a two-time-scale structure and also provides an estimate of the state derivative. The concurrent learning (CL) technique is employed to use the recorded and current data simultaneously for parameter adaptation. We show that a verifiable condition on the linear independence of the recorded data, which implies a relaxed PE condition, is sufficient to guarantee the convergence of the parameter estimation. As far as the authors' knowledge goes, the proposed CL-AEO is the first observer which is capable of jointly estimating  system state and parameter with a relaxed PE condition.
  \item An adaptive observation-based RL scheme is established for approximate optimal control of uncertain systems. The proposed CL-AEO provides the system state and model information to implement a model-based RL algorithm. Specifically, the estimated model is leveraged to evaluate the BE not only along the system trajectory, but also at any unexplored interested data points. Convergence of the developed policy to a neighborhood of the optimal policy is proved via  Lyapunov-based stability analysis. Compared with the state-of-the-art model-based RL designs \cite{Kama-2016a,Kama-2016b,Dep-2020,Dep-2020b}, our established scheme is output-feedback and does not require the state derivative information or integral calculation.
\end{enumerate}

The rest of this paper is organized as follows. Section \ref{Sec2} states the problem formulation. Section \ref{Sec3} presents the design and analysis of the CL-AEO. Section \ref{Sec4} gives the CL-AEO based RL scheme. Simulation results are provided in Section \ref{Sec5} to illustrate the effectiveness of the proposed observer and control scheme. Finally, Section \ref{Sec6} concludes the paper.

\section{Problem Formulation}\label{Sec2}

\subsection{Notations and Definitions}

Throughout the paper, big $O$-notation in terms of $\nu$ is denoted as $O(\nu)$ and it is assumed that this holds for $\nu$ positive and sufficiently small. For any continuously differentiable function $f: \mathbb{R}^{n}\times \mathbb{R}^m\rightarrow \mathbb{R}^l$, $f_x:\mathbb{R}^{n}\times \mathbb{R}^m\rightarrow \mathbb{R}^{l\times n} $ represents its gradient with respect to the first vector argument, i.e., $f_x(\nu_1, \nu_2)=\partial f(\nu_1,\nu_2)/\partial \nu_1$.  $\lambda_{\max}(P)$ and $\lambda_{\min}(P)$ denote the maximum and minimum eigenvalues of matrix $P$, respectively. $I$ denotes the identity matrix with appropriate dimension.   $\textbf{1}_A(\nu)$ is the indictor function defined by $\textbf{1}_A(\nu)=\left\{\begin{matrix} 1 ~~ \textrm{if} ~ \nu\in A, \\ 0 ~~ \textrm{if} ~ \nu \notin A. \end{matrix}\right.$ In the paper, for notation simplicity, the time variable $t$ of a signal will be omitted except when the dependence of the signal on $t$ is crucial for presentation.

The definitions of excitation and persistence of excitation of a bounded vector signal $\nu(t)$ are given as follows \cite{Chow-2013,Tao-2003}:

\emph{Definition 1:} A bounded vector signal $\nu(t)$ is \emph{exciting over an interval} $[t, t+T]$, $T>0$, $t\geq t_0$, if there exists  $\alpha>0$ such that
\begin{equation*}\label{eq3}
  \int_{t}^{t+T}\nu(\tau)\nu^{\rm{T}}(\tau)\textrm{d}\tau \geq \alpha I.
\end{equation*}

\emph{Definition 2:} A bounded vector signal $\nu(t)$ is \emph{persistently exciting (PE)} if for all $t\geq t_0$ there exist $T>0$ and  $\alpha>0$ such that
\begin{equation*}\label{eq4}
  \int_{t}^{t+T}\nu(\tau)\nu^{\rm{T}}(\tau)\textrm{d}\tau \geq \alpha I.
\end{equation*}

\subsection{Problem Statement}

Consider an $n$-dimensional nonlinear time-varying dynamic system with single-input $u$ and single-output $y$,
\begin{equation}\label{eq1}
y^{(n)}=f(y,\dot{y},\ldots,y^{(n-1)})+g(y,\dot{y},\ldots,y^{(n-1)})u,
\end{equation}
where $f(\cdot)$ and $g(\cdot)$ are continuously differentiable functions. System (\ref{eq1}) represents a wide class of physical plants, such as the wing rock phenomenon in \cite{Mona-1996} and the noncircular turning process in \cite{Wu-2009}. In this paper, we consider the case that the drift dynamics $f(\cdot)$ is in a parametric form, i.e., system (\ref{eq1}) can be written into
\begin{equation}\label{eq2}
 \left\{
  \begin{aligned}
\dot{x}= & Ax+B\left[W^{\rm{T}}\Phi(x)+g(x)u\right], \\
y= & Cx,
  \end{aligned}
\right.
\end{equation}
where $x=[x_1, \ldots, x_n]^{\rm{T}}\in\mathbb{R}^n$ is the state, $\Phi(x): \mathbb{R}^n\rightarrow \mathbb{R}^m$ is the regressor function, $W\in\mathbb{R}^m$ is the unknown constant ideal weight vector, and matrices $A\in\mathbb{R}^{n\times n}$,  $B\in\mathbb{R}^{n\times 1}$, and $C\in\mathbb{R}^{1\times n}$  represent a chain of integrators as in \cite{Khalil-2008}. To guarantee the controllability of the system, it is assumed that $g(x)$ is bounded away from zero for all $x\in\mathbb{R}^n$.
The first problem to be solved in this paper is stated as follows:

\emph{Problem 1 (Joint State-Parameter Estimation):} Given the uncertain system (\ref{eq2}), design an adaptive observer to jointly estimate the system state $x$ and parameter $W$.

The problem of joint estimation of missing state and parameter has motivated a lot of work, especially the so-called adaptive observers. However, in the existing adaptive observers such as  \cite{Farza-2018,Tyu-2013,Ibr-2018}, the restrictive PE condition is needed. The CL technique is a promising approach to relax the PE condition. But additional state derivative information is generally required \cite{Chow-2013}.  In this paper, we aim to solve Problem 1 by designing a CL-based adaptive observer, with a relaxed verifiable PE condition and without the state derivative information.

Based on the joint state-parameter estimation, we further consider the online infinite horizon optimal control problem for system (\ref{eq2}), i.e., to find the optimal control policy $u^*:\mathbb{R}^n\rightarrow \mathbb{R}$ such that the following cost functional is minimized:
\begin{equation}\label{eq6}
  J(x,u)=\int_{0}^{\infty}\left(Q(x)+u^{\rm{T}}Ru\right),
\end{equation}
where $R>0$ and $Q:\mathbb{R}^n\rightarrow \mathbb{R}^+$ is positive definite. It is well-known that this optimal control problem can be transformed into solving the following HJB equation \cite{Kirk-2004}:
\begin{align}\label{eq18}
  {V^*_x}(x)\left[Ax+B(W^{\rm{T}}\Phi(x)+g(x)u^*(x))\right] \qquad\nonumber \\
    +Q(x)+{u^*}^{\textrm{T}}(x)Ru^*(x)=0,
\end{align}
where $V^*: \mathbb{R}^n\rightarrow \mathbb{R}_{\geq 0}$, $V^*(0)=0$, is the optimal value function. The optimal control policy can be calculated from $V^*(x)$ as
\begin{equation}\label{eq9}
  u^*(x)=-\frac{1}{2}R^{-1}g^{\rm{T}}(x)B^{\rm{T}}{V_x^*}^{\rm{T}}(x).
\end{equation}
Generally speaking, an analytical solution to the HJB equation (\ref{eq18}) is not feasible, and one needs to seek an approximated solution. The second problem to be solved in this paper is then stated as follows:

\emph{Problem 2 (RL-Based Approximate Optimal Control):} Given the uncertain system (\ref{eq2}) and the cost functional (\ref{eq6}), develop a joint state-parameter estimation-based RL scheme to online approximate the optimal control policy.

Since the joint state-parameter estimation provides the system model information, the RL algorithm to be developed in this paper falls into the model-based RL community. The main advantage of the model-based RL is that by leveraging the system model, the BE can be evaluated at any points in the domain of operation, which intuitively removes the sufficient excitation assumption required by non-model-based RL algorithms.  In the following, Problems 1 and 2 will be solved in Sections \ref{Sec3} and \ref{Sec4}, respectively.

\section{Concurrent Learning-Based Adaptive Observation}\label{Sec3}

\subsection{Observer Design}

For the joint state-parameter estimation problem, similar to \cite{Farza-2018,Tyu-2013,Ibr-2018}, all the signals in system (\ref{eq2}) are assumed to be bounded. Specifically, let $x(t)\in \mathcal{X}$, $\forall t\geq 0$, where $\mathcal{X}\subset \mathbb{R}^n$ is a compact set. Before stating our candidate observer, we introduce some notations. Let $\Gamma_1=\left[\frac{l_1}{\varepsilon}, \frac{l_2}{\varepsilon^2}, \ldots, \frac{l_n}{\varepsilon^n}\right]^{\rm{T}}$ and $\Gamma_2=\frac{l_{n+1}}{\varepsilon^{n+1}}$, where $\varepsilon<1$ a small positive constant and $L=[l_1, l_2,\ldots,l_{n+1}]^{\rm{T}}\in\mathbb{R}^{n+1}$ is selected such that the following matrix is Hurwitz:
\begin{equation*}
  E=\left[
      \begin{array}{ccccc}
        -l_1 & 1 & 0 & \cdots & 0 \\
        -l_2 & 0 & 1 & \cdots & 0 \\
        \vdots &\vdots &\vdots & \ddots & \vdots \\
        -l_n & 0 & 0 & \ldots & 1 \\
        -l_{n+1} & 0 & 0 & \cdots & 0 \\
      \end{array}
    \right]\in\mathbb{R}^{(n+1)\times(n+1)}.
\end{equation*}
Let  $\varrho: \mathbb{R}\rightarrow \mathbb{R}$ be an odd smooth saturation-like function, which is characterized by $0<\varrho'(\nu)\leq 1$, $\varrho(\nu)=\nu$ if $|\nu|\leq 1$, and $\lim_{\nu\rightarrow \infty}\varrho(\nu)=1+\iota$ with $0<\iota \ll 1$ \cite{Khalil-2008}.

Let us now state the proposed concurrent learning adaptive extended observer (CL-AEO):
\begin{equation}\label{eq7}
 \left\{
  \begin{aligned}
 \dot{\widehat{x}}=& A\widehat{x}+\Gamma_1 (x_1-\widehat{x}_1)+B\left[\widehat{x}_{n+1}+g(\widehat{x})u\right],   \\
 \dot{\widehat{x}}_{n+1}=& \Gamma_2(x_1-\widehat{x}_1), \\
 \overline{x}_i= & M_i\varrho(\widehat{x}_i/M_i), ~1\leq i\leq n+1, \\
 \dot{\widehat{W}}=& \Gamma_3 \Phi(\overline{x})\left(\overline{x}_{n+1}-\widehat{W}^{\rm{T}}\Phi(\overline{x})\right) \\
 & +\sum_{j=1}^{p}\Gamma_3 \Phi(\overline{x}^{j})\left(\overline{x}^{j}_{n+1}-\widehat{W}^{\rm{T}}\Phi(\overline{x}^{j})\right),
  \end{aligned}
\right.
\end{equation}
where $\widehat{x}=[\widehat{x}_1,\ldots,\widehat{x}_{n}]^{\rm{T}}$, $\widehat{x}_{n+1}$, and $\widehat{W}$, are the estimates of $x$, $x_{n+1}\triangleq W^{\rm{T}}\Phi(x)$, and $W$, respectively; $M_i\geq \sup_{x\in \mathcal{X}}|x_i|$ are saturation bounds selected to prevent the perking phenomenon during the initial period \cite{Khalil-2008,Ran-2017a,Ran-2020}; $\overline{x}=[\overline{x}_1,\ldots, \overline{x}_n]^{\rm{T}}$ is the saturated estimate of $x$; $\Gamma_3$ is a positive definite learning rate matrix;  $j\in \{1,2,\ldots, p\}$ with $p\geq m$ denotes the index of a recorded data point, and $\overline{x}^j$ and $\overline{x}_{n+1}^j$ represent the $j$-th recorded data of $\overline{x}$ and $\overline{x}_{n+1}$, respectively.

In the sequel, several remarks are presented that provide intuitive explanations on the observer structure and reveal its properties.

\emph{Remark 1 (Two-Time-Scale Structure):} Since $\varepsilon$ is a small positive constant,  the CL-AEO (\ref{eq7}) has a two-time-scale structure in which $\widehat{x}$ and $\widehat{x}_{n+1}$ are in the fast time scale while $\widehat{W}$ is in the slow time scale. This structure is vital to implement the CL technique. Specifically, note that in the CL-AEO (\ref{eq7}), the history data of $\overline{x}$ and $\overline{x}_{n+1}$ are stored for the adaptation of $\widehat{W}$. The fast convergence of $\widehat{x}$ and $\widehat{x}_{n+1}$ guarantees the validity of the stored data. If $\widehat{x}$, $\widehat{x}_{n+1}$, and $\widehat{W}$ perform in the same time-scale, the noneligible history data will be stored. This will inevitably deteriorate the observer performance and make the convergence analysis difficult.

\emph{Remark 2 (Derivative Information):} Note that the adaptive observer (\ref{eq7}) also provides an estimate of the term $W^{\rm{T}}\Phi(x)$, which is regarded as an extended state of the system. In this case, the derivative of the $n$th state, $\dot{x}_n$, can be evaluated as $\widehat{x}_{n+1}+g(\widehat{x})u$. In the traditional CL framework  \cite{Chow-2013,Lee-2019,Kay-2019}, the state derivative information is assumed to be known or needs to be calculated by noncausal numerical smoothing techniques. These numerical processes are usually vulnerable to approximation errors, and also lay a big  barrier for the rigorous theoretical analysis \cite{Dong-2019}. In this paper, benefited from the extended design, the state derivative information is estimated simultaneously with the state. More importantly, this enables us to overcome the theoretical barrier within the traditional CL framework since the state derivative estimation error can be rigorously analyzed.

\emph{Remark 3 (Data Recording Algorithm):} In the CL-AEO (\ref{eq7}), the recorded data include the vectors $\Phi(\overline{x}^j)$ and the associated information $\overline{x}_{n+1}^j$. Let $Z=[\Phi(\overline{x}^1), \ldots, \Phi(\overline{x}^p)]$ represent the history stack, and denote $\Lambda=[\overline{x}_{n+1}^1, \ldots, \overline{x}_{n+1}^p]$. Similar to \cite{Chow-2013}, the basic idea for data recording is to update the history stack by adding data points to empty slots or by replacing an existing point if no empty slot is available to maximize the minimum singular value of $Z$. The data recording algorithm for the CL-AEO (\ref{eq7}) is given by Algorithm \ref{Alg1}.

\begin{algorithm}
 \caption{Data Recording Algorithm for the CL-AEO (\ref{eq7})}
 \label{Alg1}
 \begin{algorithmic}[1]
 \renewcommand{\algorithmicrequire}{\textbf{Input:}}
 \renewcommand{\algorithmicensure}{\textbf{Output:}}
 \STATE Set $k=1$, $Z=\textbf{0}_{m\times p}$, $\Lambda=\textbf{0}_{1\times p}$
 \IF {$k\leq p$}
  \STATE $Z(:,k)=\Phi(\overline{x})$; $\Lambda(:,k)=\overline{x}_{n+1}$
  \STATE $k=k+1$
 \ENDIF
 \IF {$k>p$}
 \STATE $Z_{\textrm{temp}}=Z$
 \STATE $S_{\textrm{old}}=\min SVD(Z^{\rm{T}})$
 \FOR {$j = 1$ to $p$}
 \STATE $Z(:,j)=\Phi(\overline{x})$
 \STATE $S(j)=\min SVD(Z^{\rm{T}})$
 \STATE $Z=Z_{\rm{temp}}$
 \ENDFOR
 \STATE find $S_{\rm{new}}=\max_{1\leq j\leq p}S(j)$ and let $j'$ denote the corresponding column index
 \IF {$S_{\rm{new}}>S_{\rm{old}}$}
 \STATE $Z(:,j')=\Phi(\overline{x})$; $\Lambda(:,j')=\overline{x}_{n+1}$
 \ENDIF
 \ENDIF
 \end{algorithmic}
 \end{algorithm}

\emph{Remark 4 (Special Cases):} In some cases, the CL-AEO (\ref{eq7}) can be constructed in more specific forms. First, if the system state $x$ is available for feedback, the structure of the CL-AEO is simplified into
\begin{equation}\label{eq-AEO1}
 \left\{
  \begin{aligned}
  \dot{\vartheta} = & \frac{l}{\varepsilon}(x_n-\vartheta)+g(x)u, ~\widehat{x}_{n+1}=\frac{l}{\varepsilon}(x_n-\vartheta), \\
 \overline{x}_{n+1}= & M_{n+1}\varrho(\widehat{x}_{n+1}/M_{n+1}),  \\
 \dot{\widehat{W}}=& \Gamma_3 \Phi(x)\left(\overline{x}_{n+1}-\widehat{W}^{\rm{T}}\Phi(x)\right) \\
 & +\sum_{j=1}^{p}\Gamma_3 \Phi(x^{j})\left(\overline{x}^{j}_{n+1}-\widehat{W}^{\rm{T}}\Phi(x^{j})\right),
  \end{aligned}
\right.
\end{equation}
where $l>0$, and $x^j$ are the $j$-th recorded state data. Second, in the case that the drift dynamics $f(x)$ is partially known, the known information can be utilized in the CL-AEO. Specifically, let $f(x)=f_0(x)+W^{\rm{T}}\Phi(x)$ with a known function $f_0$, then the $\widehat{x}$-equation in (\ref{eq7}) is modified as
\begin{equation}\label{eq-AEO2}
 \dot{\widehat{x}}= A\widehat{x}+\Gamma_1 (x_1-\widehat{x}_1)+B\left[f_0(\widehat{x})+\widehat{x}_{n+1}+g(\widehat{x})u\right].
\end{equation}

\subsection{Convergence Analysis}

The following theorem contains the convergence analysis results of the CL-AEO (\ref{eq7}).

\emph{Theorem 1:} Consider the system (\ref{eq2}) and the proposed CL-AEO (\ref{eq7}). Suppose all signals in system (\ref{eq2}) are bounded and the vector signal $\Phi(x(t))$ is exciting over a finite time interval $[0, T]$. The history stack $Z$ is empty at $t=0$, and is updated according to Algorithm \ref{Alg1} such that $\textrm{rank}(Z)=m$. Then for any $\sigma>0$ and $T_0>0$, there exists $\varepsilon^*>0$ such that $\forall \varepsilon\in(0,\varepsilon^*)$:
\begin{equation}\label{eq10}
  |x_i(t)-\widehat{x}_i(t)|\leq \sigma,  ~1\leq i\leq n+1, \forall t\geq T_0,
\end{equation}
and
\begin{equation}\label{eq11}
  \lim_{t\rightarrow \infty}\|W-\widehat{W}(t)\|\leq \sigma.
\end{equation}

\emph{Proof:} Due to the two-time-scale structure of the CL-AEO (\ref{eq7}), the proof of its convergence will be started from that of $\widehat{x}$ and $\widehat{x}_{n+1}$, and ended with the convergence of $\widehat{W}$. Consider the scaled state estimation error $\eta=[\eta_1,\ldots, \eta_{n+1}]^{\rm{T}}$ with $\eta_i=\frac{x_i-\widehat{x}_i}{\varepsilon^{n+1-i}}$, $1\leq i\leq n+1$.
 By (\ref{eq2}) and (\ref{eq7}), the dynamics of $\eta$ can be given by
\begin{equation}\label{eq12}
\dot{\eta}=\frac{1}{\varepsilon}E\eta+FW^{\textrm{T}}\Phi_x(x)\left[Ax+B\left(W^{\textrm{T}}\Phi(x)+g(x)u\right)\right],
 \end{equation}
where $F=[0 ~B^{\rm{T}}]^{\rm{T}}$. Due to the boundedness of $x$ and the continuousness of the functions $\Phi$, $\Phi_x$, and $g$,  the second term in the right-hand side of the equation above is upper bounded by an $\varepsilon$-independent positive constant $N_0$.
Let $P\in\mathbb{R}^{(n+1)\times (n+1)}$ be the unique positive definite matrix solution to the matrix equation $PE+E^{\rm{T}}P=-I$, and define a Lyapunov function candidate $V_{1}(\eta)=\eta^{\rm{T}}P\eta$. It follows that
\begin{equation}\label{eq13}
\alpha_1\|\eta\|^2 \leq  V_{1}(\eta)\leq \alpha_2\|\eta\|^2, ~\left|\frac{\partial V_1(\eta)}{\partial \eta_{n+1}}\right|\leq 2\alpha_2\|\eta\|,
\end{equation}
where $\alpha_1$ and $\alpha_2$ are the minimal and maximal eigenvalues of the matrix $P$, respectively. By some straightforward manipulations, the derivative of $V_1(\eta)$ satisfies
\begin{align}\label{eq14}
  \frac{\textrm{d}V_1(\eta)}{\textrm{d}t} \leq & -\frac{1}{\varepsilon}\|\eta\|^2+2\alpha_2N_0\|\eta\| \nonumber \\
   \leq & -\frac{1}{\alpha_2\varepsilon}V_1(\eta)+\frac{2\alpha_2N_0}{\sqrt{\alpha_1}}\sqrt{V_1(\eta)}.
\end{align}
It follows that
\begin{equation}\label{eq15}
  \frac{\textrm{d}\sqrt{V_1(\eta)}}{\textrm{d}t}\leq -\frac{1}{2\alpha_2 \varepsilon}\sqrt{V_1(\eta)}+\frac{\alpha_2N_0}{\sqrt{\alpha_1}}.
\end{equation}
By (\ref{eq13}) and (\ref{eq15}), one has
\begin{align}\label{eq16}
  \|\eta\|\leq  & \frac{\sqrt{V_1(\eta)}}{\sqrt{\alpha_1}} \nonumber \\
   \leq & \left(\frac{\sqrt{V_1(\eta(0))}}{\sqrt{\alpha_1}}-\frac{2\alpha_2^2N_0\varepsilon}{\alpha_1}\right)e^{-\frac{2\alpha_2^2N_0\varepsilon}{\sqrt{\alpha_1}}\varepsilon}+\frac{2\alpha_2^2N_0\varepsilon}{\alpha_1} \nonumber \\
   \leq & \sqrt{\frac{\alpha_2}{\alpha_1}}\|\eta(0)\|e^{-\frac{2\alpha_2^2N_0}{\sqrt{\alpha_1}\varepsilon}t}+\frac{2\alpha_2^2N_0\varepsilon}{\alpha_1}.
\end{align}
Note that the right-hand side of the inequality above is of the order of $O(\varepsilon)$ for all $t\geq t_{\varepsilon}=-(n+1)\varepsilon\ln \varepsilon$. Since $t_{\varepsilon}\rightarrow 0$ as $\varepsilon\rightarrow 0$, this proves the practical convergence of the state estimation specified by (\ref{eq10}). What is more, the saturation elements $\varrho(\widehat{x}_i/M_i)$ work in the linear zone after the convergence of the state estimation, i.e., $\overline{x}_i(t)=\widehat{x}_i(t)$, $\forall t\in[T_0,\infty)$.

Let us now consider the convergence of the parameter estimation. Denote $\widetilde{W}=W-\widehat{W}$, and define the Lyapunov function candidate $V_2(\widetilde{W})=\frac{1}{2}\widetilde{W}^{\rm{T}}\Gamma_3^{-1}\widetilde{W}$. Since $T_0$ can be made arbitrarily small and the state estimate in the parameter update law is bounded, $\widetilde{W}$ is bounded over $[0, T_0]$. Let $t_1\geq T_0$ and $t_1, \ldots, t_k, \ldots$ be a sequence where each $t_k$ denotes a time instant when the history stack $Z$ is updated. Computing the time derivative of $V_2(\widetilde{W})$ over each time interval $[t_k,t_{k+1}]$ yields
\begin{align}\label{eq17}
  \frac{\textrm{d}V_2(\widetilde{W})}{\textrm{d}t}= & -\widetilde{W}^{\rm{T}} \Phi(\overline{x})\left(\overline{x}_{n+1}-\widehat{W}^{\rm{T}}\Phi(\overline{x})\right) \nonumber \\
  & -\widetilde{W}^{\rm{T}}\sum_{j=1}^{p}\Phi(\overline{x}^{j})\left(\overline{x}^{j}_{n+1}-\widehat{W}^{\rm{T}}\Phi(\overline{x}^{j})\right) \nonumber \\
  = & -\widetilde{W}^{\rm{T}}\left[\Phi(x)\Phi^{\rm{T}}(x)+\sum_{j=1}^{p}\Phi(x^j)\Phi^{\rm{T}}(x^j)\right]\widetilde{W} \nonumber \\
  & +\delta(\widetilde{W}, \varepsilon),
\end{align}
where
\begin{align*}
  \delta(\widetilde{W},\varepsilon)= & \widetilde{W}^{\rm{T}} \Phi(x)\left(x_{n+1}-\widehat{W}^{\rm{T}}\Phi(x)\right) \nonumber \\
  & -\widetilde{W}^{\rm{T}} \Phi(\overline{x})\left(\overline{x}_{n+1}-\widehat{W}^{\rm{T}}\Phi(\overline{x})\right) \nonumber \\
  & +\widetilde{W}^{\rm{T}}\sum_{j=1}^{p}\Phi(x^{j})\left(x^{j}_{n+1}-\widehat{W}^{\rm{T}}\Phi(x^{j})\right)  \nonumber \\
  & -\widetilde{W}^{\rm{T}}\sum_{j=1}^{p}\Phi(\overline{x}^{j})\left(\overline{x}^{j}_{n+1}-\widehat{W}^{\rm{T}}\Phi(\overline{x}^{j})\right). \end{align*}
By the convergence of the state estimate, the locally Lipschitz property of the basis function $\Phi$, and some straightforward manipulations, one has $\delta(\widetilde{W},\varepsilon)\leq \iota_1\varepsilon\|\widetilde{W}\|^2+\iota_2\varepsilon\|\widetilde{W}\|$, for some $\varepsilon$-independent positive constants $\iota_1$ and $\iota_2$. This together with the fact that $\Phi(x)\Phi^{\rm{T}}(x)\geq 0$, $\forall \Phi(x)$, one gets
\begin{equation}\label{eq52}
  \frac{\textrm{d}V_2(\widetilde{W})}{\textrm{d}t} \leq -\widetilde{W}^{\rm{T}}\Psi \widetilde{W}+\iota_1\varepsilon\|\widetilde{W}\|^2+\iota_2\varepsilon\|\widetilde{W}\|,
\end{equation}
where $\Psi=\sum_{j=1}^{p}\Phi(x^j)\Phi^{\rm{T}}(x^j)\geq 0$. Note that (\ref{eq52}) guarantees  that $\widetilde{W}$ is bounded over every finite time interval $[t_k, t_{k+1}]$ if $t_{k+1}\leq T$. Recall that $\Phi(x)$ is exciting over the time interval $[0, T]$, and Algorithm \ref{Alg1} makes that the history stack $Z$ contains at least $m$ linearly independent elements for all $t\geq T$ (i.e., $\textrm{rank}(Z)=m$). Therefore for $t\geq T$, one has
\begin{equation}\label{eq53}
  \frac{\textrm{d}V_2(\widetilde{W})}{\textrm{d}t} \leq -(\lambda_{\min}\left(\Psi)-\iota_1\varepsilon\right)\|\widetilde{W}\|^2+\iota_2\varepsilon\|\widetilde{W}\|.
\end{equation}
Let $\varepsilon\in(0, \lambda_{\min}(\Psi)/\iota_1)$. Since $\lambda_{\min}(\Psi)$ is monotonically increasing, $V_2(\widetilde{W})$ is a common Lyapunov function, and consequently (\ref{eq53}) establishes the practical convergence of $\widetilde{W}$ specified by (\ref{eq11}). This completes the proof of Theorem 1. \IEEEQED

\emph{Remark 5 (Relaxed PE Condition for Parameter Estimation):} Note that from Theorem 1, a sufficient condition for guaranteing parameter estimation convergence is that  $\textrm{rank}(Z)=m$, i.e, the history stack contains $m$ linearly independent data points. This condition can be guaranteed if the system is exciting over the finite time interval when the data was recorded \cite{Chow-2013}. Compared with the existing adaptive observers \cite{Farza-2018,Tyu-2013,Ibr-2018} which require the system to be exciting over all finite intervals, this condition is much more relaxed. What is more, the rank condition only concerns with the past data and can be easily monitored online. As far as the authors' knowledge goes, the adaptive observer (\ref{eq7}) is the first attempt that addresses the joint state-parameter estimation problem with a relaxed PE condition.

\section{Reinforcement Learning-Based Approximate Optimal Control}\label{Sec4}

In this section, based on the implementation of the CL-AEO,  a simulation of experience-based RL algorithm is developed to online approximate the optimal control policy.

\subsection{Control Design}

Since the analytical solution to the HJB equation (\ref{eq18}) is generally unavailable, similar to \cite{Luo-2014,Jiang-2014}, the actor-critic NN approach is utilized to approximate the value function and the optimal control policy. According to the Weierstrass approximation theorem \cite{Courant-2008}, a continuous function can be represented by an infinite-dimensional linearly independent basis function set. In practice, one can approximate the function in a compact set with a finite-dimensional function set. Let $x\in\mathcal{X}$ and $\psi=[\psi_1(x),\ldots, \psi_r(x)]^{\rm{T}}$ be the linearly independent  continuously differentiable basis function for the value function, where $\psi_i: \mathcal{X}\rightarrow \mathbb{R}$, $1\leq i\leq r$, with $r$ the number of the neurons. Then for any given constant $\bar{\varsigma}>0$, the value function can be represented by
\begin{equation}\label{eq19}
  V^*(x)=\Theta^{\textrm{T}}\psi(x)+\varsigma(x),
\end{equation}
where $\Theta\in\mathbb{R}^r$ is the ideal weight vector, and $\varsigma:\mathcal{X}\rightarrow \mathbb{R}$ denotes the approximation error satisfying  $\sup_{x\in\mathcal{X}}|\varsigma(x)|\leq \bar{\varsigma}$ and $\sup_{x\in\mathcal{X}}|\varsigma_x(x)|\leq \bar{\varsigma}$.  Consequently, the NN representation of the idea optimal control policy is given by
\begin{equation}\label{eq20}
  u^*(x)=-\frac{1}{2}R^{-1}g^{\rm{T}}(x)B^{\rm{T}}\left(\psi_x^{\textrm{T}}(x)\Theta+\varsigma_x^{\rm{T}}(x)\right).
\end{equation}

Based on the NN representations of the value function and the optimal control policy, and using the state estimate provided by the CL-AEO, the NN-based approximation of the value function and optimal control policy are given by
\begin{align}
\label{eq21}   \widehat{V}\left(\overline{x},\widehat{\Theta}_c\right)= & \widehat{\Theta}^{\textrm{T}}_c\psi(\overline{x}), \\
\label{eq22}  \widehat{u}\left(\overline{x},\widehat{\Theta}_a\right)= & -\frac{1}{2}R^{-1}g^{\rm{T}}(\overline{x})B^{\rm{T}}\psi_x^{\rm{T}}(\overline{x})\widehat{\Theta}_a,
\end{align}
where $\widehat{\Theta}_c, \widehat{\Theta}_a\in\mathbb{R}^r$ are the weights for the critic and actor NNs, respectively.

In an RL-based controller, the main task is to design the updated laws for the NN weights by leveraging the BE as a performance metric. Traditionally, the BE is evaluated along the system trajectory, which naturally leads to a sufficient exploration requirement of the system state. In this paper, bearing in mind that the CL-AEO not only provides an estimate of the system state $x$ but also the drift dynamics $f(x)$ (i.e., $\widehat{W}^{\rm{T}}\Phi(\overline{x})$), the idea of simulation of experience-based RL \cite{Kama-2016a,Kama-2016b,Dep-2020,Dep-2020b} is employed, in which the BE is evaluated along the system trajectory, and simultaneously extrapolated to a predefined set of points $\Xi=\{x_0^i\in\mathbb{R}^n| i=1, \ldots, N\}$.  Specifically, by (\ref{eq21}) and (\ref{eq22}), the instantaneous BE evaluated along the system trajectory is given by
\begin{align}\label{eq23}
\delta_t(t) = &\widehat{V}_x(\overline{x},\widehat{\Theta}_c)\left[A\overline{x}+B\left(\widehat{W}^{\rm{T}}\Phi(\overline{x})+g(\overline{x})\widehat{u}(\overline{x},\widehat{\Theta}_a)\right)\right]\nonumber \\
  & +Q(\overline{x})+\widehat{u}^{\textrm{T}}(\overline{x},\widehat{\Theta}_a)R\widehat{u}(\overline{x},\widehat{\Theta}_a).
\end{align}
The extrapolated BE at point $x_0^i$ is given by
\begin{align}\label{eq24}
\delta_i =
&\widehat{V}_x(x_0^i,\widehat{\Theta}_c)\left[Ax_0^i+B\left(\widehat{W}^{\rm{T}}\Phi(x_0^i)+g(x_0^i)\widehat{u}(x_0^i,\widehat{\Theta}_a)\right)\right]\nonumber \\
  & +Q(x_0^i)+\widehat{u}^{\textrm{T}}(x_0^i,\widehat{\Theta}_a)R\widehat{u}(x_0^i,\widehat{\Theta}_a).
\end{align}

Then the actor and critic NNs update their weights using the BEs $\delta_t$ and $\delta_i$. A least-square update law for the critic NN is designed as
\begin{align}
\label{eq25}  \dot{\widehat{\Theta}}_c=& -k_{c1}\Gamma\frac{\mu}{\rho}\delta_t-k_{c2}\frac{\Gamma}{N} \sum_{i=1}^{N}\frac{\mu_i}{\rho_i}\delta_i,
\end{align}
where $\rho =1+\gamma\mu^{\rm{T}}\Gamma\mu$, $\rho_i =1+\gamma\mu_i^{\rm{T}}\Gamma\mu_i$, $k_{c1}, k_{c2}, \gamma>0$ are constant learning gains, and
\begin{align*}
  \mu = & \phi_x(\overline{x})\left[A\overline{x}+B\left(\widehat{W}^{\textrm{T}}\Phi(\overline{x})+g(\overline{x})\widehat{u}(\overline{x},\widehat{\Theta}_a)\right)\right], \\
  \mu_i = & \phi_x(x^i)\left[Ax_0^i+B\left(\widehat{W}^{\textrm{T}}\Phi(x_0^i)+g(x_0^i)\widehat{u}(x_0^i,\widehat{\Theta}_a)\right)\right].
\end{align*}
In (\ref{eq25}), $\Gamma:\mathbb{R}_{\geq 0}\rightarrow \mathbb{R}^{r\times r}$ represents the time-varying least-square gain matrix updated by
\begin{equation}\label{eq26}
\dot{\Gamma}=\left(\beta\Gamma-k_{c1}\frac{\Gamma\mu\mu^{\rm{T}}\Gamma}{\rho^2}\right)\textbf{1}_{\{\|\Gamma\|\leq \overline{\gamma}\}}, ~\|\Gamma(0)\|\leq \overline{\gamma},
\end{equation}
where $\beta>0$ is a constant forgetting factor and $\overline{\gamma}>0$ is a saturation constant. According to Corollary 4.3.2 in \cite{Ioan-1996}, the updated law (\ref{eq26}) guarantees that $\underline{\gamma} I \leq \Gamma(t) \leq \overline{\gamma}I$, $\forall t\geq 0$, where $\underline{\gamma}>0$. Motivated by the subsequent Lyapunov-based stability analysis, the actor NN update law is designed as
\begin{align}\label{eq27}
 \dot{\widehat{\Theta}}_a= & -k_{a1}\left(\widehat{\Theta}_a-\widehat{\Theta}_c\right)-k_{a2}\widehat{\Theta}_a
 +\frac{k_{c1}G_{t}^{\rm{T}}\widehat{\Theta}_a\mu^{\textrm{T}}}{4\rho}\widehat{\Theta}_c \nonumber \\
  & + \sum_{i=1}^{N}\frac{k_{c2}G^{\rm{T}}_{i}\widehat{\Theta}_a\mu^{\rm{T}}_i}{4N\rho_i}\widehat{\Theta}_c,
\end{align}
where $k_{a1},k_{a2}>0$ are learning gains, and
\begin{align*}
  G_{t}\triangleq & \psi_x(\overline{x})Bg(\overline{x})R^{-1}g^{\rm{T}}(\overline{x})B^{\rm{T}}\psi^{\rm{T}}_x(\overline{x}),\\
  G_{i}\triangleq & \psi_x(x_0^i)Bg(x^i)R^{-1}g^{\rm{T}}(x^i)B^{\rm{T}}\psi_x^{\rm{T}}(x_0^i).
\end{align*}

\begin{figure}
  \centering
  \includegraphics[width=0.48\textwidth,bb=10 10 520 230, clip]{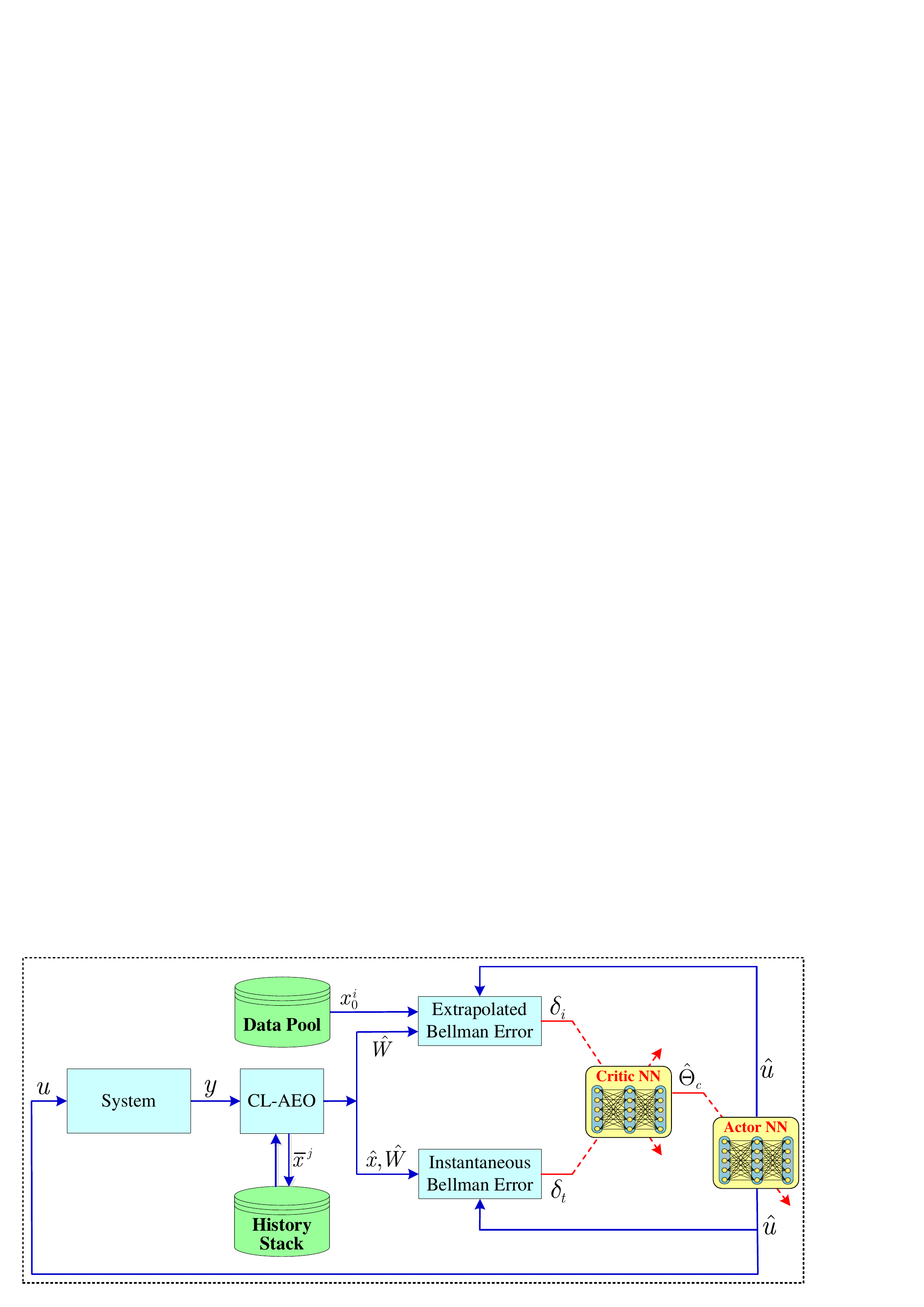}
  \caption{Adaptive observation-based efficient RL scheme.}\label{Fig1}
\end{figure}

The block diagram of the developed adaptive observation-based efficient RL scheme is illustrated in Fig. \ref{Fig1}.

\subsection{Convergence Analysis}

To facilitate the convergence analysis, the preselected data set $\Lambda$ needs to satisfy the following condition.

\emph{Assumption A1:} There exists a  constant $c>0$ such that the data points in $\Lambda$ satisfy
\begin{equation}\label{eq30}
\frac{1}{N}\inf\limits_{t\geq 0}\left(\lambda_{\min}\left\{\sum_{i=1}^N\frac{\mu_i\mu_i^{\rm{T}}}{\rho_i}\right\}\right)\geq c.
\end{equation}

\emph{Remark 6 (Relaxed PE Condition for RL-Based Control):} Note that unlike the standard PE condition, the excitation condition (\ref{eq30}) can be monitored online. What is more, by leveraging the estimated system model, the BE can be extrapolated to any selected data point. Therefore, the excitation condition (\ref{eq30}) can be met heuristically by selecting more data points than the number of neurons, i.e., $N \gg r$ \cite{Kama-2016a}. In practice, to fulfill Assumption A1, the data points $x_0^i$, $1\leq i\leq N$, can be select on an $\underbrace{a\times a \cdots \times a}_{n}$ data grid which covers the interested domain of operation $\mathcal{X}\subseteq\mathbb{R}^n$, where $a$ is an appropriately large positive integer. \IEEEQED

For subsequent use, let us specify the compact set $\mathcal{X}=\{x\in\mathbb{R}^n; \|x\|\leq \tau_x+1\}$ and define $\mathcal{X}_0=\{x\in\mathbb{R}^n; \|x\|\leq \tau_x\}$, where $\tau_x>0$. Define a concatenated state $Z(t)=[x^{\rm{T}}(t), \widetilde{W}^{\rm{T}}(t), \widetilde{\Theta}_c^{\rm{T}}(t),  \widetilde{\Theta}_a^{\rm{T}}(t)]^{\rm{T}}$, and functions $V_c=\frac{1}{2}\widetilde{\Theta}_c^{\rm{T}}\Gamma^{-1}\widetilde{\Theta}_c$ and $V_a=\frac{1}{2}\widetilde{\Theta}_a^{\rm{T}}\widetilde{\Theta}_a$.  Denote $\tau_{W}=\frac{1}{2}\Gamma_3^{-1}\|\widetilde{W}(0)\|^2$, $\tau_c=\frac{1}{2}\underline{\gamma}^{-1}\|\widetilde{\Theta}_c(0)\|^2+1$, $\tau_a=$ $\max\{V_a(\widetilde{\Theta}_a(0)),V_a(\widetilde{\Theta}_a)_{\|\widetilde{\Theta}_a\|\geq (\iota_5+\iota_{10})/|\iota_9|}\}+1$, where $\iota_5$, $\iota_9$, and $\iota_{10}$ will be specified latter.  Define several compact sets:
\begin{align*}
\Omega_{W}^0=&\{\widetilde{W}\in\mathbb{R}^m; V_2(\widetilde{W})\leq \tau_W\}, \\
\Omega_{W}^1=&\{\widetilde{W}\in\mathbb{R}^m; V_2(\widetilde{W})\leq \tau_W+1\}, \\
\Omega_{c}^0=&\{\widetilde{\Theta}_c\in\mathbb{R}^r; V_c(\widetilde{\Theta}_v)\leq \tau_c\}, \\
\Omega_{c}^1=&\{\widetilde{\Theta}_c\in\mathbb{R}^r; V_c(\widetilde{\Theta}_v)\leq \tau_c+1\}, \\
\Omega_{a}^0=&\{\widetilde{\Theta}_a\in\mathbb{R}^r; V_a(\widetilde{\Theta}_a)\leq \tau_a\}, \\
\Omega_{a}^1=&\{\widetilde{\Theta}_a\in\mathbb{R}^r; V_a(\widetilde{\Theta}_a)\leq \tau_a+1\}, \\
\Omega^0= & \mathcal{X}_0\times \Omega^0_W\times \Omega^0_c\times \Omega^0_a, \\
\Omega^1= & \mathcal{X}\times\Omega^1_W \times \Omega^1_c\times \Omega^1_a.
\end{align*}

\emph{Theorem 2:} Consider the closed-loop system formed by plant (\ref{eq2}), CL-AEO (\ref{eq7}), control (\ref{eq22}), and RL update laws (\ref{eq25})-(\ref{eq27}). The history stack satisfies $\textrm{rank}(Z)=m$, and is updated according to Algorithm \ref{Alg1}. Suppose Assumption A1 is satisfied and $x(0)\in \mathcal{X}_0-\partial \mathcal{X}_0$. Then there exists $\varepsilon^\dag>0$ such that for any $\varepsilon\in(0,\varepsilon^\dag)$:
\begin{itemize}
  \item the CL-AEO is convergent in the sense of (\ref{eq10}) and (\ref{eq11});
  \item the state $x$ and the NN weight estimation errors $\widetilde{\Theta}_c\triangleq \widehat{\Theta}_c-\Theta$ and $\widetilde{\Theta}_a\triangleq \widehat{\Theta}_a-\Theta$ are uniformly ultimately bounded.
\end{itemize}

\emph{Proof:} Since $x(0)\in \mathcal{X}_0-\partial \mathcal{X}_0$, $\widetilde{W}(0)\in \Omega_W^0-\partial \Omega_W^0$, $\widetilde{\Theta}_c(0)\in \Omega_c^0-\partial \Omega_c^0$, and $\widetilde{\Theta}_a(0)\in \Omega_a^0-\partial \Omega_a^0$, the initial concatenated state $Z(0)$ is an interior point of $\Omega^0$. This together with the fact that the state estimate output of the CL-AEO (\ref{eq7}) is bounded yields that there exists an $\varepsilon$-independent $t_0>0$ such that $Z(t)\in\Omega^0$, $\forall t\in[0,t_0]$. Since $\Omega^0\subseteq \Omega^1-\partial \Omega^1$, we let $Z(t)\in \Omega^1$, $\forall t\in[0,t_1]$, where $t_1>t_0$. In the following we will show that for sufficiently small $\varepsilon$, $t_1$ can be selected as infinity.

In the time interval $[0,t_1]$, since $Z(t)$ is bounded, one can follow a same line as the arguments as in (\ref{eq12})-(\ref{eq16}) to conclude that for any $t_0'\in (t_0, t_1)$, there exists $\varepsilon_1>0$ such that for any $\varepsilon\in(0,\varepsilon_1)$, $\|\eta(t)\|=O(\varepsilon)$, $\forall t\in[t_0',t_1]$. It follows that $\overline{x}_i(t)=\widehat{x}_i(t)$, $1\leq i\leq n+1$, $\forall t\in[t_0',t_1]$.

To facilitate the subsequent analysis, the NN-based approximations $\widehat{V}_x(\overline{x},\widehat{\Theta}_c)$ and $\widehat{u}(\overline{x},\widehat{\Theta}_a)$ are expressed as
\begin{align}
  \widehat{V}_x\left(\overline{x},\widehat{\Theta}_c\right)= & \Theta^{\rm{T}} \psi_x(\widehat{x})+\left(\widehat{\Theta}_c-\Theta\right)^{\rm{T}}\psi_x(\widehat{x}) \nonumber \\
\label{eq31}  =& V^*_x(\widehat{x})+\widetilde{\Theta}^{\rm{T}}_c\psi_x(\widehat{x})-\varsigma_x(\widehat{x}),\\
    \widehat{u}\left(\overline{x},\widehat{\Theta}_a\right)
  =& -\frac{1}{2}R^{-1}g^{\rm{T}}(\widehat{x})B^{\rm{T}}\psi^{\rm{T}}_x(\widehat{x})\left(\Theta+\widetilde{\Theta}_a\right) \nonumber\\
  =& u^*(\widehat{x})-\frac{1}{2}R^{-1}g^{\rm{T}}(\widehat{x})B^{T}\psi_x^{\rm{T}}(\widehat{x})\widetilde{\Theta}_a \nonumber\\
\label{eq32a}   &  +\frac{1}{2}R^{-1}g^{\rm{T}}(\widehat{x})B^{\rm{T}}\varsigma^{\rm{T}}_x(\widehat{x}).
\end{align}
By inserting (\ref{eq19}), (\ref{eq20}), (\ref{eq31}), and (\ref{eq32a}) into (\ref{eq23}),  the instantaneous BE $\delta_t$ can be written as
\begin{align}\label{eq28}
\delta_t = & \mu^{\rm{T}}\widetilde{\Theta}_c-V_x^*(\widehat{x})B\left[\widetilde{W}^{\rm{T}}\Phi(\widehat{x})+g(\widehat{x})(u^*(\widehat{x})-\widehat{u}(\widehat{x},\widehat{\Theta}_a))\right] \nonumber \\
  & -\varsigma_x(\widehat{x})\left[A\widehat{x}+B\left(\widehat{W}\Phi(\widehat{x})+g(\widehat{x})\widehat{u}(\widehat{x},\widehat{\Theta}_a)\right)\right] \nonumber \\
  & +u^*(\widehat{x})\left[-g^{\rm{T}}(\widehat{x})B^{T}\psi_x^{\rm{T}}(\widehat{x})\widetilde{\Theta}_a+g^{\rm{T}}(\widehat{x})B^{\rm{T}}\varsigma^{\rm{T}}_x(\widehat{x})\right] \nonumber \\
  & +\frac{1}{4R}\left[-g^{\rm{T}}(\widehat{x})B^{T}\psi_x^{\rm{T}}(\widehat{x})\widetilde{\Theta}_c+g^{\rm{T}}(\widehat{x})B^{\rm{T}}\varsigma^{\rm{T}}_x(\widehat{x})\right]^2\nonumber \\
= & \mu^{\rm{T}}\widetilde{\Theta}_c-\Theta^{\rm{T}}\psi_x(\widehat{x})B\Phi(\widehat{x})\widetilde{W} \nonumber \\
& +\frac{1}{2}\Theta^{\rm{T}}\psi_x(\widehat{x})Bg(\widehat{x})R^{-1}g^{\rm{T}}(\widehat{x})B^{\rm{T}}\varsigma^{\rm{T}}_x(\widehat{x}) \nonumber \\
& -\varsigma(\widehat{x})\left(A\widehat{x}+BW^{\rm{T}}\Phi(\widehat{x})\right) \nonumber \\
  & +\frac{1}{4}\widetilde{\Theta}_a^{\rm{T}}\psi_x(\widehat{x})Bg(\widehat{x})R^{-1}g^{\rm{T}}(\widehat{x})B^{T}\psi_x^{\rm{T}}(\widehat{x})\widetilde{\Theta}_a\nonumber \\
  & +\frac{1}{4}\varsigma_x(\widehat{x})Bg(\widehat{x})R^{-1}g^{\rm{T}}(\widehat{x})B^{\rm{T}}\varsigma^{\rm{T}}_x(\widehat{x}) \nonumber \\
 \triangleq &  \mu^{\rm{T}}\widetilde{\Theta}_c-\Theta^{\rm{T}}\psi_x(\widehat{x})B\Phi(\widehat{x})\widetilde{W}+\frac{1}{4}\widetilde{\Theta}^{\rm{T}}_aG_{t}\widetilde{\Theta}_a+\Delta_t,
\end{align}
where
\begin{align*}
 \Delta_t= & \frac{1}{2}\Theta^{\rm{T}}\psi_x(\widehat{x})Bg(\widehat{x})R^{-1}g^{\rm{T}}(\widehat{x})B^{\rm{T}}\varsigma^{\rm{T}}_x(\widehat{x}) \nonumber \\
 & -\varsigma(\widehat{x})\left(A\widehat{x}+BW^{\rm{T}}\Phi(\widehat{x})\right) \nonumber \\
 & +\frac{1}{4}\varsigma_x(\widehat{x})Bg(\widehat{x})R^{-1}g^{\rm{T}}(\widehat{x})B^{\rm{T}}\varsigma^{\rm{T}}_x(\widehat{x}).
\end{align*}
Similarly, the extrapolated BE evaluated at the selected point $x_0^i$ can be expressed as
\begin{equation}\label{eq33}
\delta_i=\mu_i^{\rm{T}}\widetilde{\Theta}_c-\Theta^{\rm{T}}\psi_x(x_0^i)B\Phi(x_0^i)\widetilde{W}+\frac{1}{4}\widetilde{\Theta}^{\rm{T}}_aG_{i}\widetilde{\Theta}_a+\Delta_i,
\end{equation}
where
\begin{align*}
 \Delta_i= & \frac{1}{2}\Theta^{\rm{T}}\psi_x(x_0^i)Bg(x_0^i)R^{-1}g^{\rm{T}}(x_0^i)B^{\rm{T}}\varsigma^{\rm{T}}_x(x_0^i) \nonumber \\
 & -\varsigma(x_0^i)\left(Ax_0^i+BW^{\rm{T}}\Phi(x_0^i)\right) \nonumber \\
 & +\frac{1}{4}\varsigma_x(x_0^i)Bg(x_0^i)R^{-1}g^{\rm{T}}(x_0^i)B^{\rm{T}}\varsigma^{\rm{T}}_x(x_0^i).
\end{align*}

Let us consider the Lyapunov function candidate given by
\begin{equation}\label{eq32}
  \mathcal{V}(Z)=V^*(x)+V_{2}(\widetilde{W})+V_c(\widetilde{\Theta}_c)+V_a(\widetilde{\Theta}_a).
\end{equation}
The derivative of $V^*(x)$ can be computed as
\begin{align}\label{eq43}
  \dot{V}^*(x)= & V_x^*(x)[Ax+B(W^{\rm{T}}\Phi(x)+g(x)\widehat{u}(\widehat{x},\widehat{\Theta}_a))] \nonumber \\
  = & V_x^*(x)\left[Ax+B\left(W^{\rm{T}}\Phi(x)+g(x)u^*(x)\right)\right]\nonumber \\
    & +V_x^*(x)B\left[g(x)\left(\widehat{u}(\widehat{x},\widehat{\Theta}_a)-u^*(\widehat{x})\right) \right.\nonumber \\
    & +g(x)\left(u^*(\widehat{x})-u^*(x)\right)\Big].
\end{align}
By (\ref{eq18}) and (\ref{eq32a}), the locally Lipschitz property of the functions $V_x^*$, $\psi_x$, $g$, and $u^*$, and the facts that $\|x(t)-\widehat{x}(t)\|=O(\varepsilon^2)$ and $x(t), \widehat{x}(t)\in \mathcal{X}$, $\forall t\in[t_0',t_1]$,  the derivative of $V^*(x)$ satisfies
\begin{equation}\label{eq44}
\dot{V}^*(x)\leq -\iota_3\|x\|^2+\iota_4\varepsilon\|\eta\|+\iota_5\|\widetilde{\Theta}_a\|+\iota_6\overline{\varsigma},
\end{equation}
where $\iota_3=\lambda_{\min}(Q)$, and $\iota_4$ to $\iota_6$ are $\varepsilon$-independent positive constants.

By the RL update laws (\ref{eq25})-(\ref{eq27}), the derivative of $V_c(\widetilde{\Theta}_c)+V_a(\widetilde{\Theta}_a)$ satisfies
\begin{align}\label{eq32}
    & \dot{V}_c(\widetilde{\Theta}_c)+\dot{V}_a(\widetilde{\Theta}_a) \nonumber \\
  = & \widetilde{\Theta}_c^{\rm{T}}\left(-k_{c1}\frac{\mu}{\rho}\delta_t-\frac{k_{c2}}{N} \sum_{i=1}^{N}\frac{\mu_i}{\rho_i}\delta_i\right)\nonumber \\
  &  -\frac{1}{2}\widetilde{\Theta}_c^{\rm{T}}\Gamma^{-1}\left(\beta\Gamma-k_{c1}\frac{\Gamma\mu\mu^{\rm{T}}\Gamma}{\rho^2}\right)\Gamma^{-1}\widetilde{\Theta}_c \nonumber \\
  & +\widetilde{\Theta}_a^{\rm{T}}\left(-k_{a1}(\widehat{\Theta}_a-\widehat{\Theta}_c)-k_{a2}\widehat{\Theta}_a\right)\nonumber\\
  & +\widetilde{\Theta}_a^{\rm{T}}\left(\frac{k_{c1}G_{t}^{\rm{T}}\widehat{\Theta}_a\mu^{\textrm{T}}}{4\rho} + \sum_{i=1}^{N}\frac{k_{c2}G^{\rm{T}}_{i}\widehat{\Theta}_a\mu^{\rm{T}}_i}{4N\rho_i}\right)\widehat{\Theta}_c.
  & \end{align}
Inserting  (\ref{eq28}) and (\ref{eq33}) into (\ref{eq32}) leads to
\begin{align}\label{eq45}
  & \dot{V}_c(\widetilde{\Theta}_c)+\dot{V}_a(\widetilde{\Theta}_a) \nonumber \\
  = & -k_{c1}\widetilde{\Theta}_c^{\rm{T}}\frac{\mu\mu^{\rm{T}}}{\rho}\widetilde{\Theta}_c +k_{c1}\widetilde{\Theta}_c^{\rm{T}}\frac{\mu}{\rho}\Theta^{\rm{T}}\psi_x(\widehat{x})B\Phi(\widehat{x})\widetilde{W} \nonumber  \\
   & -k_{c1}\widetilde{\Theta}_c^{\rm{T}}\frac{\mu}{\rho}\Delta_t-\frac{k_{c2}}{N}\widetilde{\Theta}_c^{\rm{T}}\sum_{i=1}^N\frac{\mu_i\mu_i^{\rm{T}}}{\rho_i}\widetilde{\Theta}_c -\widetilde{\Theta}_c^{\rm{T}}\frac{k_{c2}}{N}\sum_{i=1}^N\frac{\mu_i}{\rho_i}\Delta_i \nonumber \\
   & +\frac{k_{c2}}{N}\widetilde{\Theta}_c^{\rm{T}}\sum_{i=1}^N\frac{\mu_i}{\rho_i}\Theta^{\rm{T}}\psi_x(x_0^i)B\Phi(x_0^i)\widetilde{W} \nonumber \\
   &  -\frac{1}{2}\beta\widetilde{\Theta}_c^{\rm{T}}\Gamma^{-1}\widetilde{\Theta}_c+k_{c1}\widetilde{\Theta}_c^{\rm{T}}\frac{\mu\mu^{\rm{T}}}{\rho^2}\widetilde{\Theta}_c \nonumber \\
   & -(k_{a1}+k_{a2})\|\widetilde{\Theta}_a\|^2+k_{a1}\widetilde{\Theta}_a^{\rm{T}}\widetilde{\Theta}_c-k_{a2}\widetilde{\Theta}_a^{\rm{T}}\Theta \nonumber \\
   & +\left(\|\Theta\|^2\widetilde{\Theta}_a^{\rm{T}}+\|\widetilde{\Theta}_a\|^2\Theta^{\rm{T}}+\widetilde{\Theta}_a^{\rm{T}}\Theta\widetilde{\Theta}_c\right)\nonumber \\
   & \times \left(\frac{k_{c1}G_{t}^{\rm{T}}\mu^{\textrm{T}}}{4\rho} + \sum_{i=1}^{N}\frac{k_{c2}G^{\rm{T}}_{i}\mu^{\rm{T}}_i}{4N\rho_i}\right).
\end{align}
Since $\underline{\gamma}I\leq \Gamma(t)\leq \overline{\gamma}I$, the normalized regressor $\frac{\mu}{\rho}$ is bounded as $\|\frac{\mu}{\rho}\|\leq \frac{1}{2\sqrt{\gamma \underline{\gamma}}}$, and so does $\frac{\mu_i}{\rho_i}$.  These together with Assumption A1 yields
\begin{align}\label{eq46}
       & \dot{V}_c(\widetilde{\Theta}_c)+\dot{V}_a(\widetilde{\Theta}_a) \nonumber \\
  \leq & \frac{k_{c1}}{\gamma\underline{\gamma}}\|\widetilde{\Theta}_c\|^2+\iota_7\|\widetilde{W}\|^2+\frac{k_{c1}^2}{2\iota_7\gamma\underline{\gamma}}\|\Theta\|^2\vartheta_1^2\|\widetilde{\Theta}_c\|^2  \nonumber \\
  & -k_{c2}c\|\widetilde{\Theta}_c\|^2+\frac{k_{c2}^2}{2\iota_7\gamma\underline{\gamma}}\|\Theta\|^2\vartheta_2^2\|\widetilde{\Theta}_c\|^2 \nonumber \\
  & -\frac{\beta}{2\overline{\gamma}}\|\widetilde{\Theta}_c\|^2+\frac{k_{c1}}{4\gamma\underline{\gamma}}\|\widetilde{\Theta}_c\|^2-(k_{a1}+k_{a2})\|\widetilde{\Theta}_a\|^2 \nonumber \\
  & +\frac{k_{a1}}{2}(\|\widetilde{\Theta}_c\|^2+\|\widetilde{\Theta}_a\|^2)+k_{a2}\|\Theta\|\|\widetilde{\Theta}_a\| \nonumber \\
  & +\vartheta_3\|\Theta\|^2\|\widetilde{\Theta}_a\|+\vartheta_3\|\Theta\|\|\widetilde{\Theta}_a\|^2 \nonumber \\
  & +\frac{\vartheta_3\|\Theta\|}{2}(\|\widetilde{\Theta}_c\|^2+\|\widetilde{\Theta}_a\|^2)+\frac{k_{c1}}{2\sqrt{\gamma\underline{\gamma}}}\|\widetilde{\Theta}_c\||\Delta_t| \nonumber \\
  & +\frac{k_{c2}}{2\sqrt{\gamma\underline{\gamma}}}\|\widetilde{\Theta}_c\|\max_{1\leq i\leq N}|\Delta_i| \nonumber \\
 \triangleq &  \iota_7\|\widetilde{W}\|^2+\iota_8\|\widetilde{\Theta}_c\|^2+ \iota_9\|\widetilde{\Theta}_a\|^2+\iota_{10}\|\widetilde{\Theta}_a\|+\iota_{11},
\end{align}
where
\begin{align*}
  \iota_7 \triangleq & \frac{1}{2}\lambda_{\min}(\Psi), \\
  \iota_8 \triangleq & \frac{5k_{c1}}{4\gamma\underline{\gamma}}+\frac{k_{c1}^2}{2\iota_7\gamma\underline{\gamma}}\|\Theta\|^2\vartheta_1^2-k_{c2}c+\frac{k_{c2}^2}{2\iota_7\gamma\underline{\gamma}}\|\Theta\|^2\vartheta_2^2 \\
  & -\frac{\beta}{2\overline{\gamma}}+\frac{k_{a1}}{2}+\frac{\vartheta_3\|\Theta\|}{2}, \\
  \iota_9 \triangleq & -\frac{k_{a1}}{2}-k_{a2}+\frac{3}{2}\vartheta_3\|\Theta\|, \\
  \iota_{10} \triangleq & k_{a2}\|\Theta\|+\vartheta_3\|\Theta\|^2, \\
  \iota_{11} \triangleq & \frac{k_{c1}}{2\sqrt{\gamma\underline{\gamma}}}\|\widetilde{\Theta}_c\||\Delta_t|+\frac{k_{c2}}{2\sqrt{\gamma\underline{\gamma}}}\|\widetilde{\Theta}_c\|\max_{1\leq i\leq N}|\Delta_i|, \\
  \vartheta_1\triangleq  & \sup_{\widehat{x}\in\mathcal{X}}\|\psi_x(\widehat{x})B\Phi(\widehat{x})\|, \\
  \vartheta_2\triangleq  & \max_{1\leq i\leq N}\|\psi_x(x_0^i)B\Phi(x_0^i)\|, \\
  \vartheta_3\triangleq  & \frac{k_{c1}}{8\sqrt{\gamma\underline{\gamma}}}\sup_{\widehat{x}\in\mathcal{X}}\|G_t\|+\frac{k_{c2}}{8\sqrt{\gamma\underline{\gamma}}}\max_{1\leq i\leq N}\|G_i\|.
\end{align*}

By (\ref{eq53}), (\ref{eq44}), and (\ref{eq46}), one gets
\begin{align}\label{eq47}
  \dot{\mathcal{V}}\leq & -\iota_3\|x\|^2-(\iota_7-\iota_1\varepsilon)\|\widetilde{W}\|^2+\iota_2\varepsilon\|\widetilde{W}\|+\iota_8\|\widetilde{\Theta}_c\|^2 \nonumber \\
  & +\iota_9\|\widetilde{\Theta}_a\|^2+(\iota_5+\iota_{10})\|\widetilde{\Theta}_a\|+\iota_{12},
\end{align}
where $\iota_{12}\triangleq \iota_4\varepsilon\|\eta\|+\iota_6\overline{\varsigma}+\iota_{11}$. The sufficient conditions for the learning gains are given by
\begin{align}\label{eq48}
  \beta>\frac{\overline{\gamma}}{\iota_7\gamma\underline{\gamma}}\|\Theta\|^2\left(k_{c1}^2\vartheta_1^2+k_{c2}^2\vartheta_2^2\right)+\overline{\gamma}\vartheta_3\|\Theta\|,
\end{align}
\begin{equation}\label{eq49}
  k_{c2}>\frac{5k_{c1}}{4\gamma\underline{\gamma}c}+\frac{k_{a1}}{2c},  ~k_{a2}>\frac{3}{2}\vartheta_3\|\Theta\|.
\end{equation}
Recall that 1) the system state estimate is convergent in the sense that $\|\eta(t)\|=O(\varepsilon)$, $\forall t\in[t_0', t_1]$; 2) $\Psi>0$ and $\lambda_{\min}(\Psi)$ is monotonically incerasing according to Algorithm \ref{Alg1}; 3) the upper bound of the approximation error $\overline{\varsigma}$ can be made arbitrarily small by increasing the number of NN neurons; and 4) the concatenated state $Z(t) \in \Omega^1$, $\forall t\in[0, t_1]$. Therefore, provided the learning gains satisfying (\ref{eq48}) and (\ref{eq49}), the derivative of $\mathcal{V}$ satisfies
\begin{equation}\label{eq48b}
  \dot{\mathcal{V}}(Z(t))\leq 0,  ~t\in[t'_0,t_1].
\end{equation}
This indicates that $t_1$ can be selected as infinity, and consequently the uniformly ultimately boundedness of $x$, $\widetilde{\Theta}_c$, and $\widetilde{\Theta}_v$ are guaranteed. The convergence of the CL-AEO follows a same line of the arguments as in Theorem 1, with the boundedness of the signals in the closed-loop system. Finally, (\ref{eq47}) also implies that the system state $x$ converges to the neighbourhood of origin, and the actor NN weight $\widehat{\Theta}_a$ approximates the ideal weight $\Theta$ if $(\iota_5+\iota_{10})/|\iota_9|$ is sufficiently small (i.e., $k_{a1}$ is selected sufficiently large). This completes the proof of Theorem 2. \IEEEQED

\emph{Remark 7 (Efficient RL Scheme):} In this paper, the ideas from adaptive observers and CL  are leveraged to design efficient model-based RL algorithms. The advantages of our proposed approach are threefold: 1) It adopts relaxed and verifiable PE conditions. For the joint state-parameter estimation process, the CL technique is employed to replace the PE conditions required in the existing adaptive observers \cite{Farza-2018,Tyu-2013,Ibr-2018} as a verifiable rank condition. For the approximated optimal control process, the PE condition is relaxed via the simulation of experience technique. 2) It is numerically efficient since it does not require the state derivative information required in \cite{Chow-2013,Kama-2016a,Kama-2016b,Dep-2020} and the integral calculation required in \cite{Dep-2020b}. 3) It is output feedback and does not require the probing signal. The RL algorithms  in \cite{Lewis-2009,Jiang-2012,Lewis-2010,Luo-2014,Jiang-2014,Jiang-2015,Bai-2020,Wen-2017,Ma-2019,Lewis-2014,Xue-2020,Yang-2019,Kama-2016a,Kama-2016b,Dep-2020,Dep-2020b} require full state-feedback  and in \cite{Jiang-2012,Lewis-2010,Jiang-2014,Jiang-2015,Ma-2019,Lewis-2014} require the probing signal. \IEEEQED

\section{Simulation Study}\label{Sec5}

In this section, simulations are provided to demonstrate the effectiveness and superiority of the proposed CL-AEO and learning-based controller.
We consider a second-order uncertain nonlinear system given by
\begin{equation}\label{eq-sim}
 \left\{
  \begin{aligned}
\dot{x}_1=& x_2, \\
\dot{x}_2= & -x_1+W_1x_2+W_2(x_1+x_2)(\sin(x_1)+2)^2 \\
 & +(\sin(x_1)+2)u, \\
 y = & x_1,
  \end{aligned}
\right.
\end{equation}
where $W_1=-1.5$, $W_2=0.5$. The simulations consist of two parts. In the first part, the control signal $u$ is manually selected to verify that the proposed CL-AEO removes the requirement of the PE condition.  In the second part, the closed-loop performance of the adaptive observation-based RL algorithm is illustrated.

\subsection{Adaptive Observation}

\begin{figure}[!t]
\centering
\begin{subfigure}[b]{0.42\textwidth}
 \centering
 \includegraphics[width=1.0\textwidth,bb=5 0 315 370, clip]{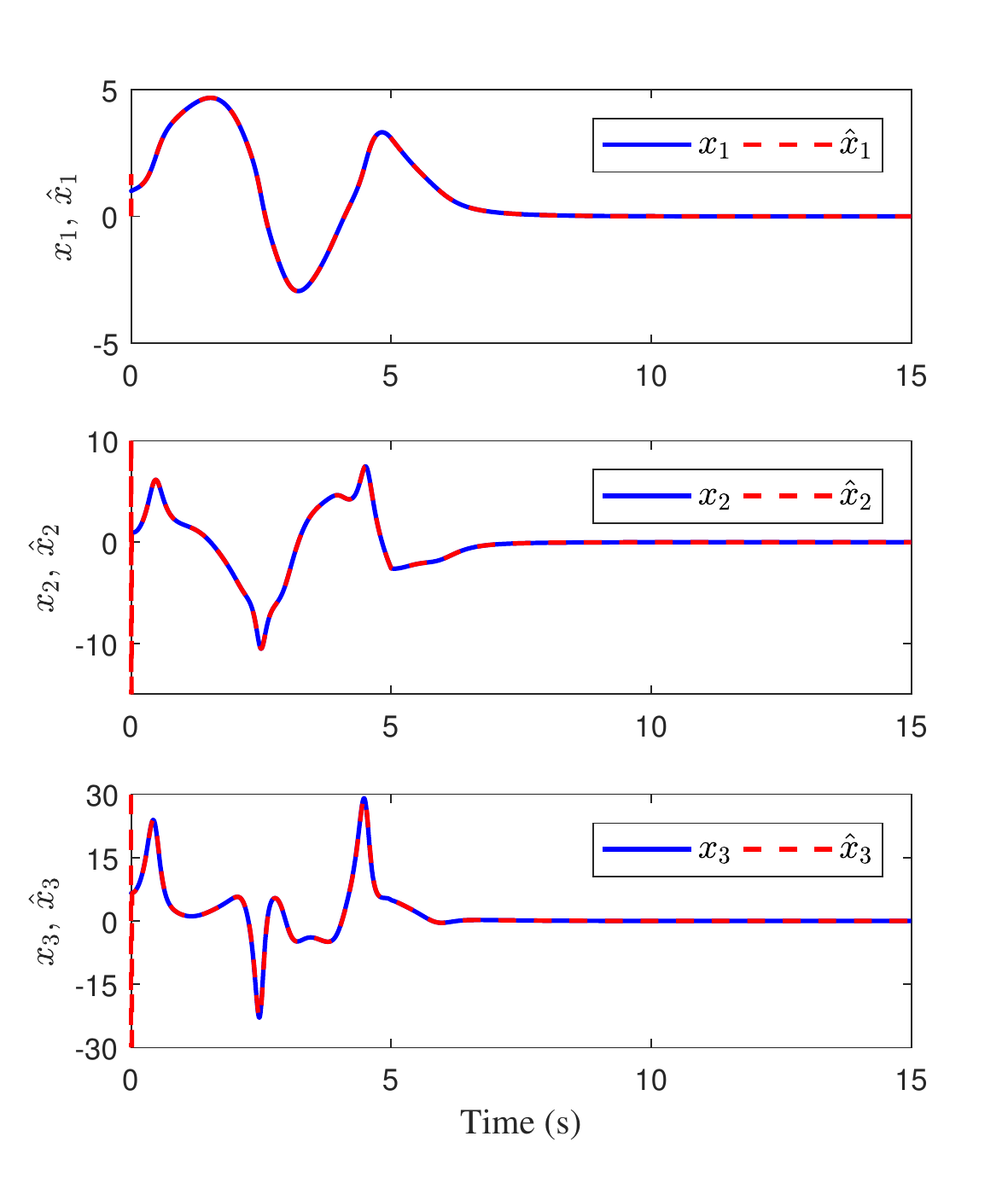}
 \caption{Estimation of the system state $x$ and extended state $x_3$.} \label{fig2a}
        \end{subfigure}
        %\hfill
        \vskip\baselineskip
        \begin{subfigure}[b]{0.42\textwidth}
 \centering
 \includegraphics[width=1.0\textwidth,bb=5 0 315 245, clip]{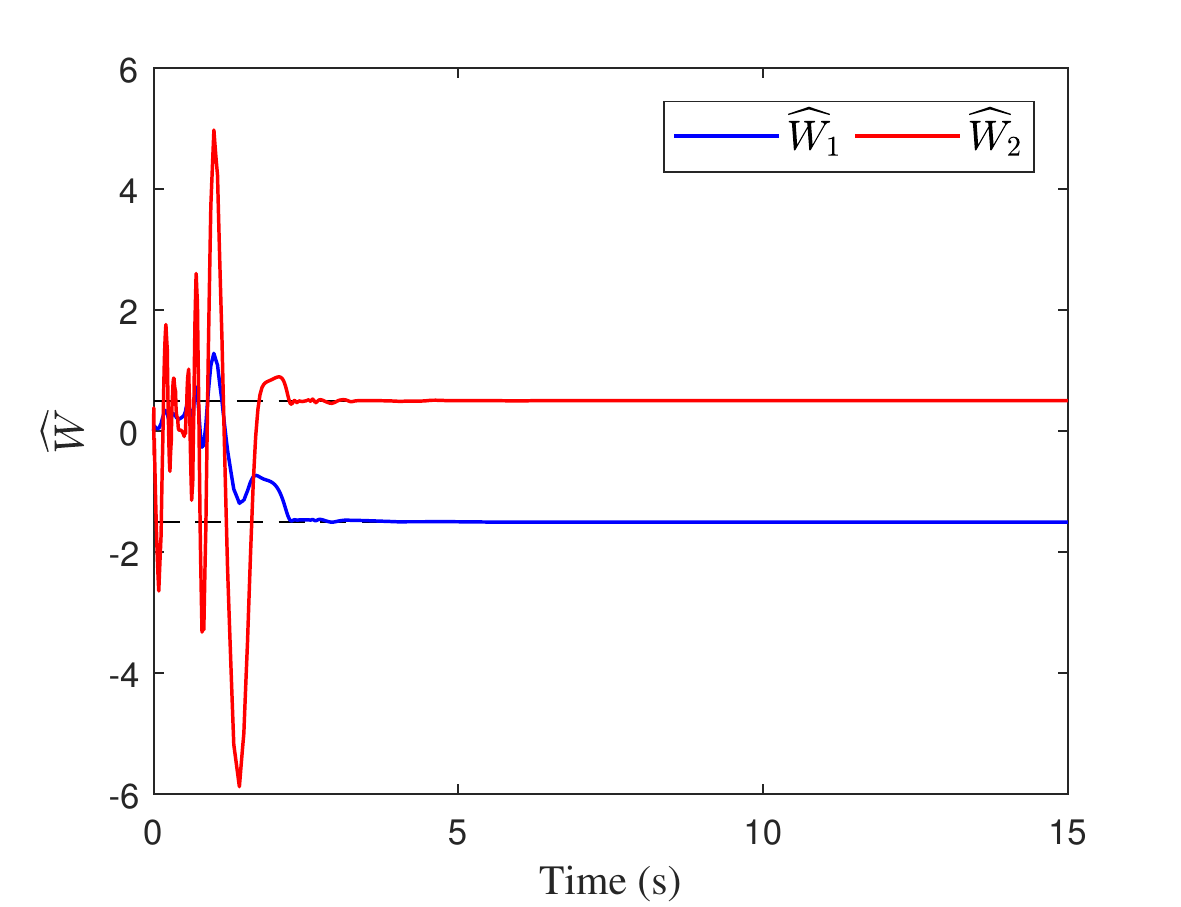}
 \caption{Estimation of the system parameter $W$. The true values of $W$ are shown as dashed lines.}\label{fig2b}
        \end{subfigure}
        \caption{Response of the CL-AEO with control signal $u_1$.}
        \label{fig2}
\end{figure}

\begin{figure}[!t]
\centering
\begin{subfigure}[b]{0.42\textwidth}
 \centering
 \includegraphics[width=1.0\textwidth,bb=3 0 315 245, clip]{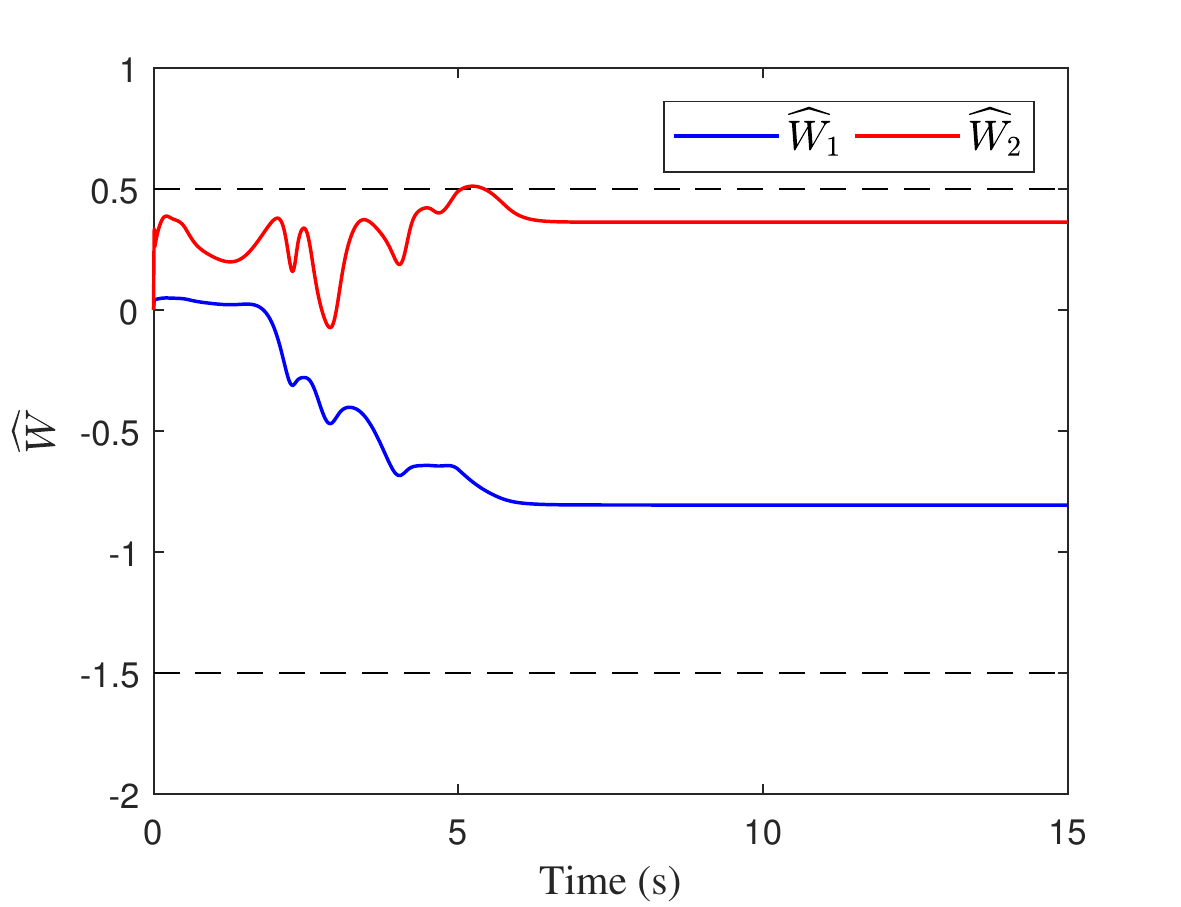}
 \caption{With control signal $u_1$.}\label{fig3a}
        \end{subfigure}
        %\hfill
        \vskip\baselineskip
        \begin{subfigure}[b]{0.42\textwidth}
 \centering
 \includegraphics[width=1.0\textwidth,bb=3 0 315 245, clip]{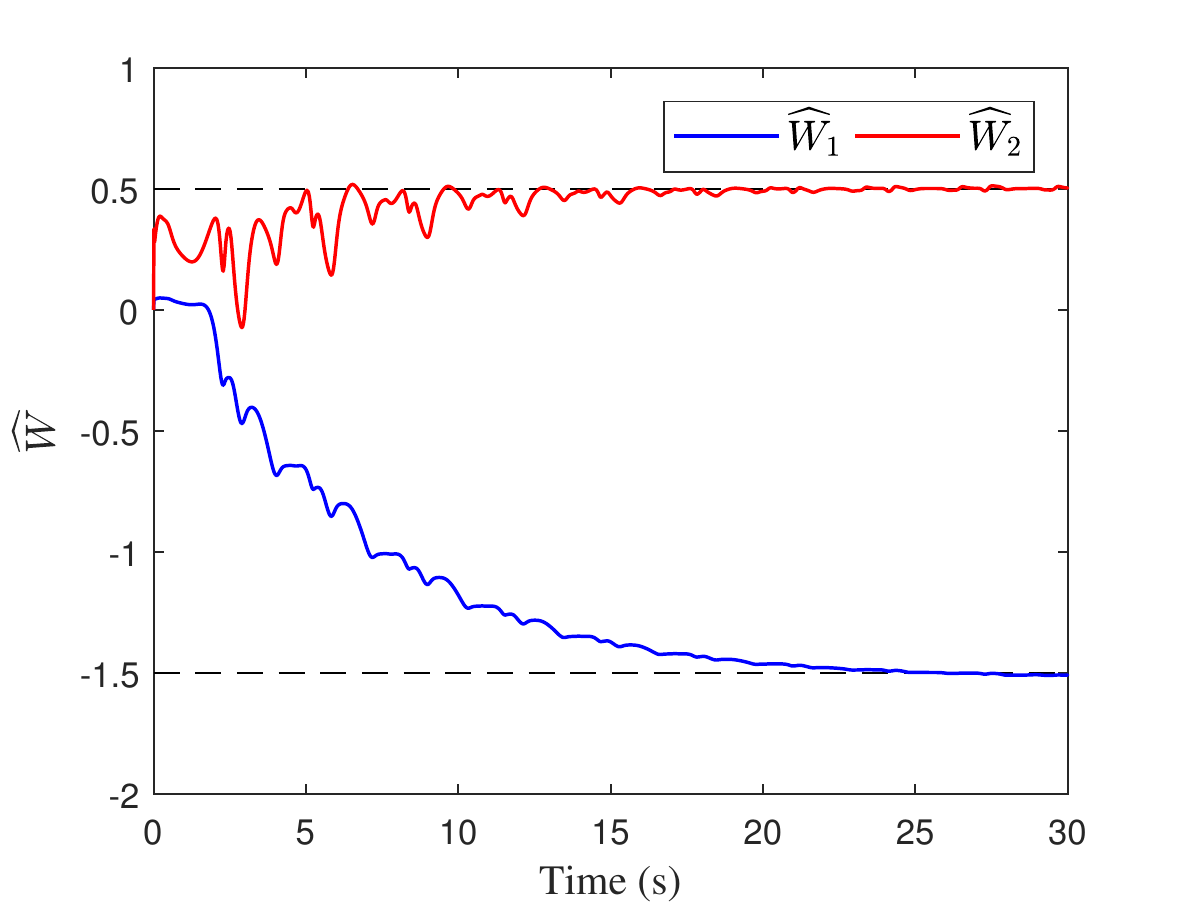}
 \caption{With control signal $u_2$.}\label{fig3b}
        \end{subfigure}
        \caption{System parameter estimations of the AEO with different control signals.}
        \label{fig3}
\end{figure}

\begin{figure}[!t]
\centering
\begin{subfigure}[b]{0.42\textwidth}
 \centering
 \includegraphics[width=1.0\textwidth,bb=5 0 315 253, clip]{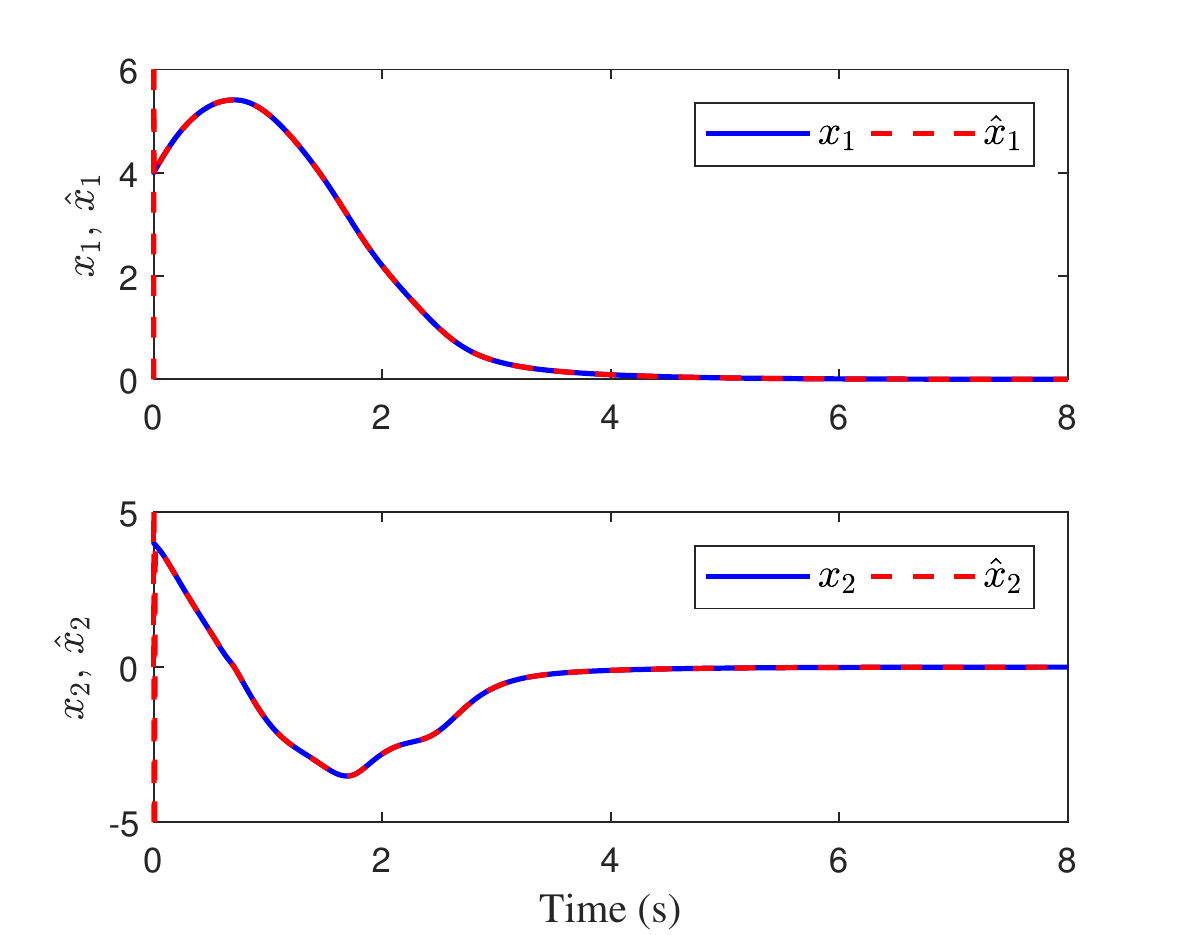}
 \caption{Trajectories of the system sate $x$ and its estimate $\widehat{x}$.} \label{fig4a}
        \end{subfigure}
        %\hfill
        \vskip\baselineskip
        \begin{subfigure}[b]{0.42\textwidth}
 \centering
 \includegraphics[width=1.0\textwidth,bb=0 0 315 245, clip]{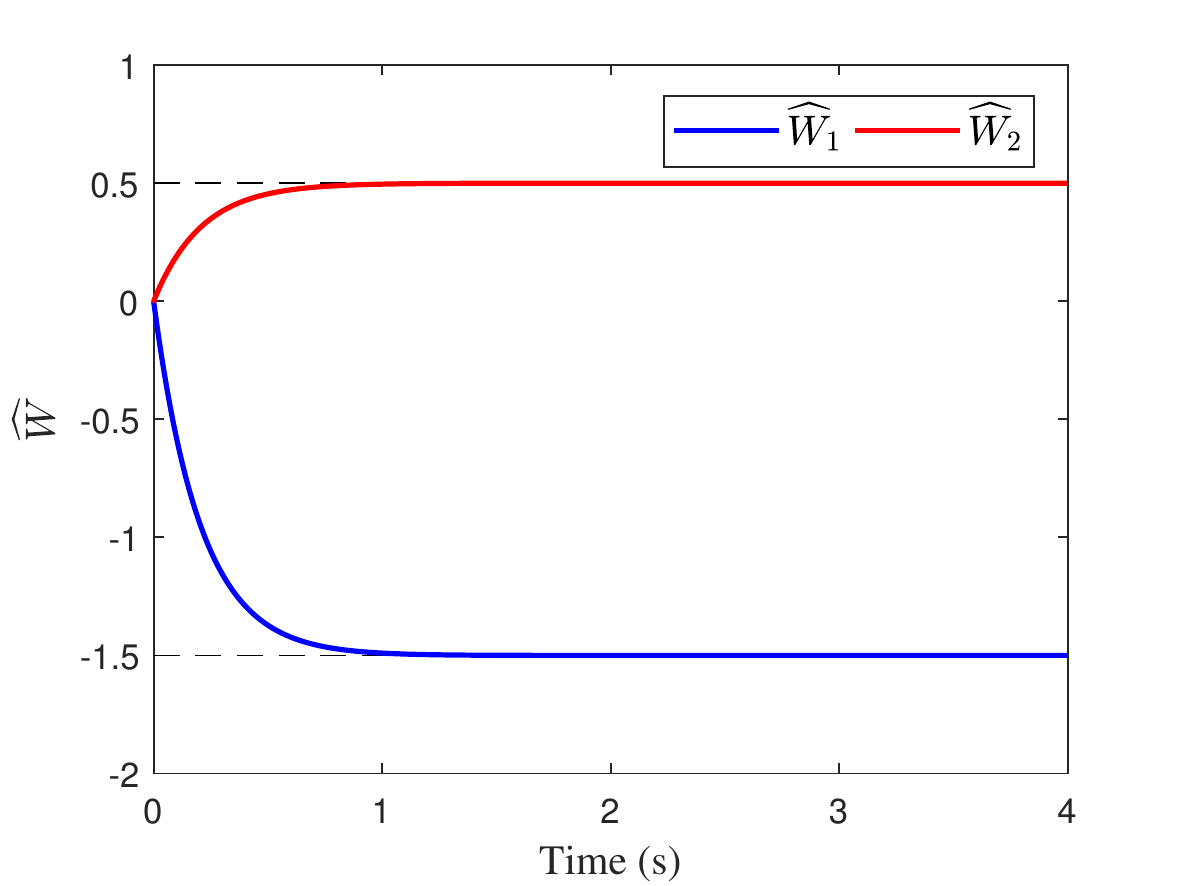}
 \caption{Trajectories of the parameter estimate $\widehat{W}$.}\label{fig4b}
        \end{subfigure}
        \caption{Responses of the system and the CL-AEO.}
        \label{fig4}
\end{figure}

In this subsection, two control signals are considered:
\begin{align}\label{eq50}
  u_1(t)=&\left\{\begin{aligned}
  & -0.9(\cos(2x_1)+2)(x_1+x_2) \\
  & \quad +10\sin(4\pi t),  ~0\leq t \leq 5,\\
  & -0.9(\cos(2x_1)+2)(x_1+x_2),  ~t>5,
  \end{aligned}\right. \\
  u_2(t)=&-0.9(\cos(2x_1)+2)(x_1+x_2) \nonumber \\
\textbf{\label{eq50}}  & +10\sin(4\pi t),  t\geq 0.
\end{align}
Note that the control signal $u_1$ only  makes the system finitely exciting, while the control signal $u_2$ persistently exciting. The system initial condition is set as $x(0)=[1 ~1]$. The CL-AEO is designed with $L=[3 ~3 ~1]$, $\varepsilon=0.001$, $M_1=5$, $M_2=12$, $M_3=30$, $p=5$, and $\Gamma_3=0.4I$. The instantaneous data for recording is selected according to Algorithm \ref{Alg1}. For comparison, the  counterpart adaptive extended observer (AEO) without concurrent learning is also simulated. The initial conditions of the observers are all set as  0.

Fig. \ref{fig2} shows the performance of the proposed CL-AEO under the control signal $u_1$. One can see that both the system state $x$, the extended state $x_3\triangleq W_1x_2+W_2(x_1+x_2)(\sin(x_1)+2)^2$, and the parameter $W$ are estimated satisfactorily. This figure also illustrates the two-time-scale  property of the proposed CL-AEO, i.e., the convergence speed of $\widehat{x}$ and $\widehat{x}_3$ is faster than $\widehat{W}$.  Figs. \ref{fig3a} and \ref{fig3b} show  the system parameter estimations of the AEO with finite exciting control signal $u_1$ and persistently exciting control signal $u_2$,  respectively. It can be seen that the observer without current learning cannot arrive at the ideal weight $W$ with a finite exciting control signal. With a persistently exciting control signal, the AEO arrives at the ideal weight $W$, it is however, with a lower convergence speed. These illustrate the advantages of our proposed CL-AEO.

% that the convergence speed of the state estimate can be sufficiently fast, while the convergence speed of the parameter estimate is limited. The reason for this distinction is that the convergence of the parameter estimate can only be achieved when the system is exciting  over a finite interval, i.e., the concurrent learning algorithm collects sufficient rich data such that the rank condition in Assumption A? is satisfied.

\subsection{Approximate Optimal Control}

To consider the optimal control problem, the cost functional given by (\ref{eq6}) is specified   with $Q=x^{\rm{T}}\overline{Q}x$ and $R=1$, where $\overline{Q}=\left[
                \begin{array}{cc}
                  2 & 1 \\
                  1 & 1 \\
                \end{array}
              \right]$
is positive definite. This cost functional is selected because the corresponding optimal control problem has a known analytical solution. According to the procedure in \cite{Nev-1996}, the optimal value function and control policy are $V^*(x)=1.5x_1^2+2x_1x_2+x_2^2$ and $u^*(x)=-(\cos(2x_1)+2)(x_1+x_2)$, respectively.

We consider the scenario that the system initial state $x(0)$ locates in the interior of the compact set $\{x\in\mathbb{R}^2; \|x\|\leq 10\}$. For the CL-AEO, the settings are the same as those in the previous subsection except $M_1=10$, $M_2=10$, and $M_3=100$. What is more, the stack $Z$ is initialized with three history data points recorded in the previous simulation.
For the RL controller, the basis function is selected as $\psi(x)=[x_1^2 ~x_1x_2 ~x_2^2]^{\rm{T}}$, which implies that the ideal weight $\Theta=[1.5 ~2 ~1]^{\rm{T}}$; the learning gains are selected as $k_{c1}=1$, $k_{c2}=5$, $k_{a1}=80$, $k_{a2}=0.1$, $\gamma=0.5$, and $\beta=100$; the upper bound of the norm of $\Gamma$ is set as $\overline{\gamma}=1000$; the data set $\Lambda$ to extrapolate the BE contains $121$ data points which located on a $11\times 11$ data grid covers the domain $[-10 ~10]\times [-10 ~10]$ (i.e., $x_1$ and $x_2$ are both selected every 0.5 from -10 to 10). Simulation is done with initial conditions $x(0)=[4 ~4]^{\rm{T}}$, $[\widehat{x}_1(0) ~\widehat{x}_2(0) ~\widehat{x}_3(0)]^{\rm{T}}=[0 ~0 ~0]^{\rm{T}}$, $\widehat{W}(0)=0$, $\widehat{\Theta}_c(0)=\widehat{\Theta}_a(0)=[0 ~0 ~0]^{\rm{T}}$, and $\Gamma(0)=\textrm{diag}\{100,100,100\}$.

\begin{figure}[!t]
 \centering
 \includegraphics[width=0.42\textwidth,bb=5 0 315 240, clip]{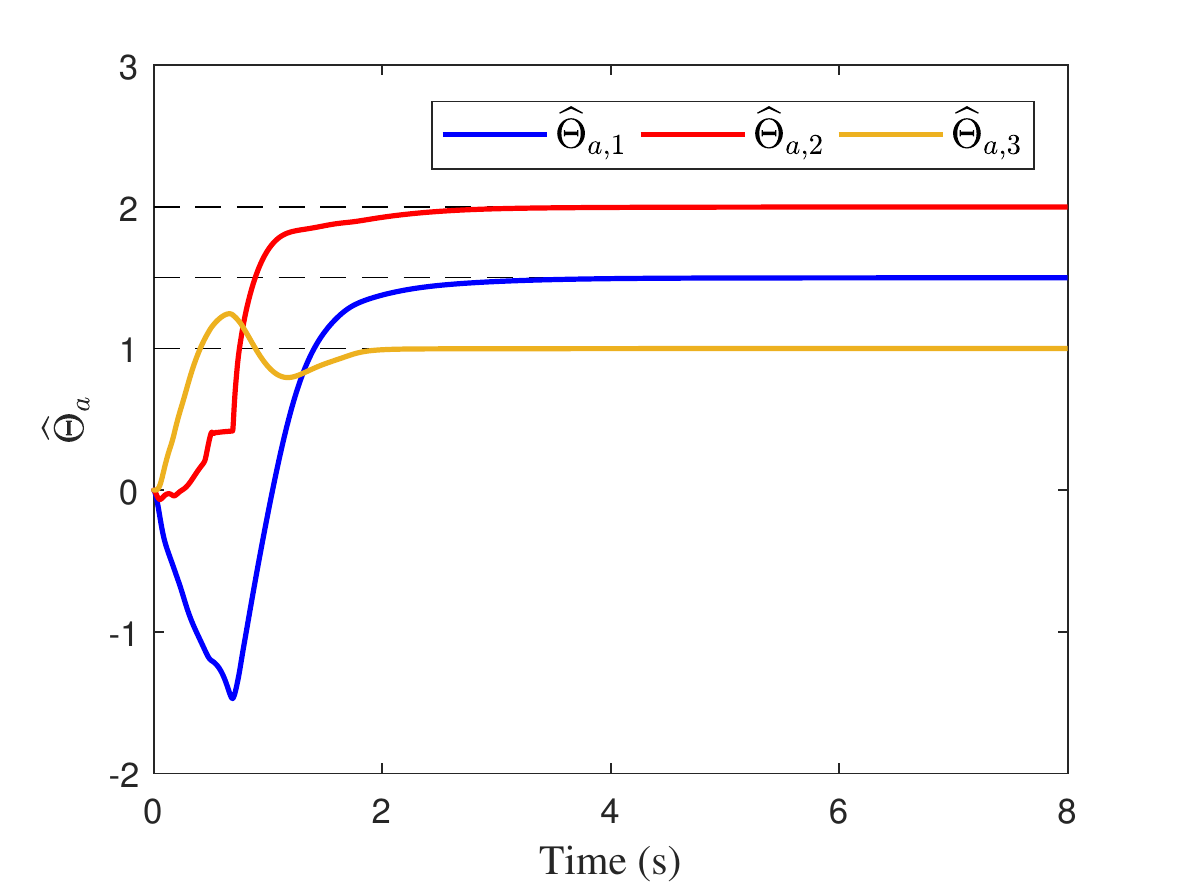}
 \caption{Trajectories of the actor NN weight $\widehat{\Theta}_a$. The true values of $\Theta_a$ are shown as dashed lines.}\label{fig5}
\end{figure}

\begin{figure}
 \centering
 \includegraphics[width=0.42\textwidth,bb=5 0 315 200, clip]{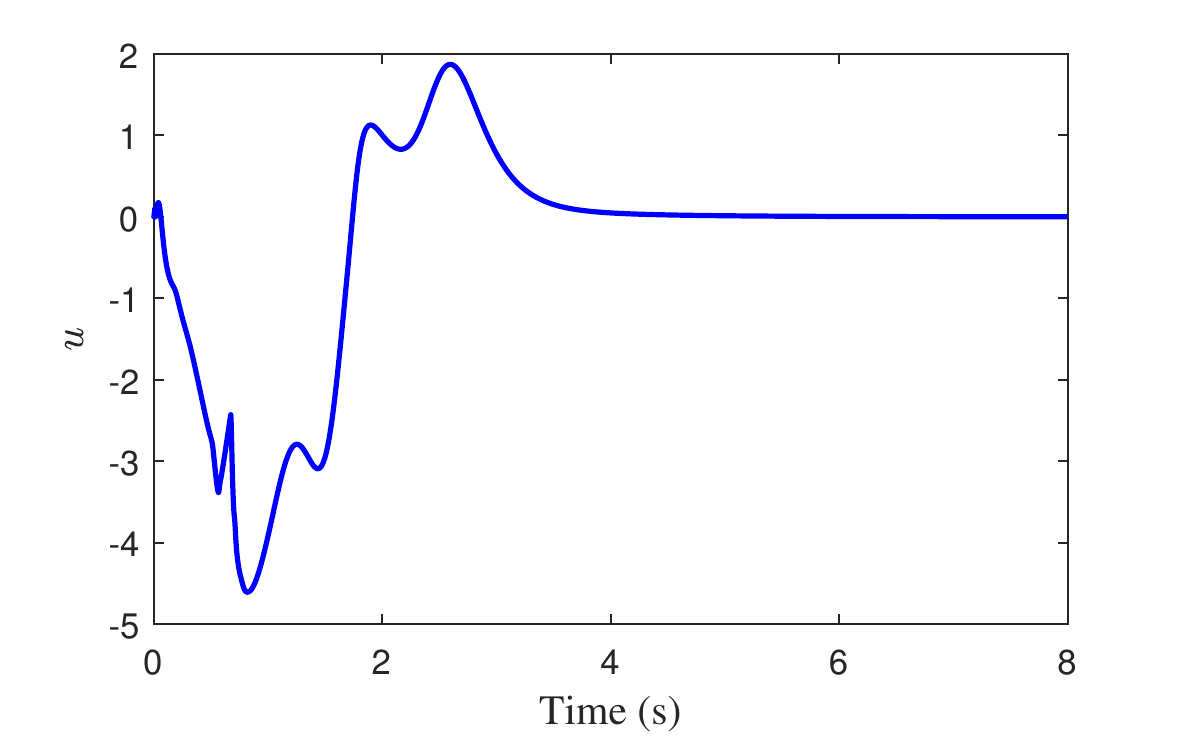}
 \caption{Trajectories of the control input.}\label{fig6}
\end{figure}

Fig. \ref{fig4} shows the responses of the system and the CL-AEO, from which one can see that the learning-based controller regulates the system state to 0, and the CL-AEO provides accurate joint estimation of the system state  and parameter. Fig. \ref{fig5} illustrates that the actor NN weight $\widehat{\Theta}_a$ converges to its real value, and this indicates that the learning-based controller performs an approximate optimal control property.
Fig. \ref{fig6} depicts the control signal. The trajectory of $\frac{1}{N}\left(\lambda_{\min}\left\{\sum_{i=1}^N\frac{\mu_i\mu_i^{\rm{T}}}{\rho_i}\right\}\right)$ is plotted in Fig. \ref{fig7}, from which one can see that the condition in Assumption A1 is satisfied.

Finally, we illustrate the advantage of the implementation of the simulation of experience mechanism. Fig. \ref{fig8} showes the trajectories of the actor NN weight $\widehat{\Theta}_a$ without simulation of experience, i.e., the BE is evaluated only along the system trajectory. From this figure, one can see that $\widehat{\Theta}_a$ cannot converge to its real value. The reason is that the  system state doesn't explore sufficient points in the state space. In this case, to make sure  $\widehat{\Theta}_a\rightarrow \Theta$, a carefully selected probing signal is required to inject into the system \cite{Jiang-2012,Lewis-2010,Jiang-2014,Jiang-2015,Ma-2019,Lewis-2014}, which will deteriorate the closed-loop performance. In our simulation of experience-based RL, by leveraging the estimated system model, the BE can be extrapolated to any selected data point, and hence the probing signal is not required anymore.

\section{Conclusion}\label{Sec6}

An adaptive observation-based efficient RL framework is established for uncertain systems, which consists of a CL-AEO to provide the system state and model information, and a simulation-of-experience-based RL mechanism to approximate the optimal control policy. Both the observation and control processes adopt relaxed and verifiable PE conditions. The obtained results provide a novel and practical adaptive observation-based solution for the implementation of model-based RL, since it is output feedback and does not require the state derivative information or integral calculation..

In practice, the system control coefficient function $g(x)$ maybe also uncertain. Therefore, for future works, we aim at extending the developed approach to systems with uncertain control coefficient functions.

\end{document}